\documentclass[12pt]{amsart} 
\usepackage{fullpage}
\usepackage{mathtools}
\usepackage[all]{xy}
\usepackage[inline]{asymptote}
\usepackage{palatino, mathpazo}
\usepackage{newtxtext,newtxmath}
\usepackage{hyperref}

\usepackage{enumerate}
\usepackage{graphicx}
\usepackage{bbm}

\newcommand{\GL}{\mathrm{GL}}

\newcommand{\DDD}{\mathcal{D}}
\newcommand{\MMM}{\mathcal{M}}
\newcommand{\OOO}{\mathcal{O}}

\newcommand{\LLL}{\mathcal{L}}

\newcommand{\FF}{\mathbb{F}}
\newcommand{\CC}{\mathbb{C}}

\newcommand{\GG}{\mathbb{G}}
\newcommand{\HH}{\mathbb{H}}
\newcommand{\DD}{\mathbb{D}}
\newcommand{\RR}{\mathbb{R}}

\newcommand{\ZZ}{\mathbb{Z}}

\newcommand{\PP}{\mathbb{P}}
\newcommand{\QQ}{\mathbb{Q}}

\newcommand{\Ppp}{\mathfrak{p}}

\newcommand{\tto}[1]{\xrightarrow{#1}}
\newcommand{\gen}[1]{\left\langle#1\right\rangle}

\newcommand{\pq}[1]{\left(#1\right)}
\newcommand{\aq}[1]{\left|#1\right|}

\newcommand{\cq}[1]{\left\{#1\right\}}
\newcommand{\bq}[1]{\left[#1\right]}

\newcommand{\ol}[1]{\overline{#1}}
\newcommand{\del}{\partial}

\newcommand{\FFF}{\mathcal{F}}

\newcommand{\inv}{^{-1}}
\newcommand{\f}[2]{\frac{#1}{#2}}

\newcommand{\sm}{\setminus}

\newcommand{\du}{^{*}}
\DeclareMathOperator{\Gal}{Gal}

\DeclareMathOperator{\Aut}{Aut}
\DeclareMathOperator{\Ext}{Ext}

\DeclareMathOperator{\Spec}{Spec}

\newcommand{\CLie}{\mathcal{CL}}

\newcommand{\simplecap}[2]{\refstepcounter{figure}\vskip1emFigure \ref{#1}. #2\label{#1}}%
\newcommand{\Cor}{\text{\rm Cor}}

\newcommand{\pg}[2]{\big(#1\;|\;#2\big)}

\newcommand{\ldbr}{\left[\hspace{-0.02in}\left[\,}
\newcommand{\rdbr}{\,\right]\hspace{-0.02in}\right]}

\newcommand{\HHH}{\mathcal{H}}

\newcommand{\nil}{\text{\rm nil}}
\newcommand{\Hod}{\text{\rm Hod}}

\newcommand{\MT}{\text{\rm MT}}
\newcommand{\gr}{\text{\rm gr}}

\newcommand{\Lie}{\text{\rm Lie}}
\newcommand{\Der}{\text{\rm Der}}
\DeclareMathOperator{\Sym}{Sym}
\newcommand{\reg}{\text{\rm reg}}

\newcommand{\Mot}{\text{\rm Mot}}
\newcommand{\End}{\text{\rm End}}

\newcommand{\PM}{\mathcal{PM}}

\newcommand{\HS}{\mathrm{HS}}

\newcommand{\MHS}{\mathrm{MHS}}

\def\o{\omega}
\def\ob{\ol\omega}

\newtheorem{thm}{Theorem}
\newtheorem{lma}[thm]{Lemma}
\newtheorem{cly}[thm]{Corollary}
\newtheorem{conj}[thm]{Conjecture}

\newcommand{\CL}{\mathcal{CL}}

\newcommand{\uns}{^\times}
\newcommand{\chain}{^\bullet}
\renewcommand{\div}{{\rm div}}
\newcommand{\av}{{\rm av}}
\def\ch{\chain}

\newcommand{\MM}{\mathcal{MM}}

\title[Motivic $\pi_1$ of CM elliptic curves and geometry of Bianchi hyperbolic threefolds]{Motivic fundamental groups of CM elliptic curves\\ and geometry of Bianchi hyperbolic threefolds}
\author{Nikolay Malkin}

\begin{document}

\begin{abstract}
    In this paper we describe a connection between realizations of the action of the motivic Galois group on the motivic fundamental groups of Gaussian and Eisenstein elliptic curves punctured at the $\Ppp$-torsion points, $\pi_1^{\Mot}(E-E[\Ppp],v_0)$, and the geometry of the Bianchi hyperbolic threefolds $\Gamma_1(\Ppp)\sm\HH^3$, where $\Gamma_1(\Ppp)$ is a congruence subgroup of $\GL_2(\End(E))$. The first instance of such a connection was found by A.Goncharov in \cite{goncharov-euler}.

    In particular, we study the Hodge realization of the image of the above action in the fundamental Lie algebra, a pronilpotent Lie algebra carrying a filtration by \emph{depth}. The depth-1 associated graded quotient of the image is fully described by Beilinson and Levin's elliptic polylogarithms (\cite{beilinson-levin}). In this paper, we consider the depth-2 associated graded quotient. One of our main results is the construction of a homomorphism from the complex computing the cohomology of a certain local system on the Bianchi threefold to this quotient's standard cochain complex. This result generalizes those of \cite{goncharov-euler}, as well the connection between modular manifolds and the motivic fundamental group of $\GG_m$ punctured at roots of unity (\cite{goncharov-polylogs-modular,goncharov-motivic-modular}).
    
    Our construction uses the mechanism of Hodge correlators, canonical generators in the image that were defined in \cite{goncharov-hodge-correlators}. Our second main result is a system of \emph{double shuffle relations} on the canonical real periods of these generators, the Hodge correlator integrals. These relations deform the relations on the depth-2 Hodge correlators on the projective line, previously found by the author in \cite{malkin-shuffle}.
\end{abstract}

\maketitle

\setcounter{tocdepth}{2}
\tableofcontents

\section{Introduction}

\label{sec:intro}

In this paper we describe a connection between the realizations of motivic fundamental groups of CM elliptic curves and the geometry of Bianchi hyperbolic threefolds. The first instance of this connection was described by A.Goncharov in \cite{goncharov-euler}.

\subsubsection*{Motivation} We aim to study the action of the motivic Galois group on the motivic fundamental group of an elliptic curve punctured at the $\Ppp$-torsion points with tangential base point $v_0$:
\begin{equation}
    \Gal_\Mot\circlearrowright\pi_1^\Mot(E-E[\Ppp],v_0).
    \label{eqn:gal_mot_action}
\end{equation}
The objects in (\ref{eqn:gal_mot_action}) are still conjectural, but we can study them in their realizations. The results of this paper are in the Hodge realization. However, the picture is easiest to introduce in the $\ell$-adic realization.

As a running example, take $E$ to be the CM elliptic curve $E=\CC/(\ZZ+\ZZ[i])$ and $\Ppp\subset\ZZ[i]$ an ideal. The $\ell$-adic realization of the motivic fundamental group, $\pi_1^{(\ell)}(E-E[\Ppp],v_0)$, is simply the pro-$\ell$ completion of the topological fundamental group $\pi_1(E-E[\Ppp],0)$. It is equipped with an action of the absolute Galois group $\Gal(\ol\QQ/\QQ)$ by automorphisms.

The Maltsev construction (\cite{deligne-p1p1}, \S9) makes out of $\pi_1^{(\ell)}(E-E[\Ppp],v_0)$ a pro-$\ell$ Lie algebra $A_{E,\Ppp}$ over $\QQ_\ell$, generated by $H_1(E;\ZZ)$ and loops around the punctures in $E[\Ppp]$. It carries two filtrations: by \emph{weight} and by \emph{depth}. The increasing weight filtration $W$ (see \cite{deligne-hodge1}) is invariant under the Galois action, and the geometric Frobenius element %
acts on $\gr^wA_{E,\Ppp}$ with eigenvalues of norm $\ell^{w/2}$. The descreasing depth filtration $D$ is defined by the lower central series of the linearization of \[\ker\pq{\pi_1^{(\ell)}(E-E[\Ppp],v_0)\to\pi_1^{(\ell)}(E,v_0)}.\]

The filtrations $W$ and $D$ induce filtrations on $\End(A_{E,\Ppp})$, and, by restriction, on the image of the action of $\Gal(\ol\QQ/\QQ)$. Taking its associated graded Lie algebra for the weight filtration, we obtain a graded Lie algebra $\Lie_{(\ell)}(E,E[\Ppp])$, the \emph{elliptic Galois Lie algebra}. We study the quotient of this Lie algebra induced by the quotient of $A_{E,\Ppp}$ by the adjoint action of $H_1(E;\ZZ)$, and take the coinvariants of the translation action of $E[\Ppp]$ on $E$ (amounting to averaging the base point).%
This quotient is called the \emph{symmetric} Galois Lie algebra $\Lie_{(\ell)}^{\rm sym}(E,E[\Ppp])$.

The structure of the depth-$d$ graded quotients of $\Lie_{(\ell)}^{\rm sym}(E,E[\Ppp])$ is well understood in depths 0 and 1. In depth 0, this algebra simply vanishes. In depth 1, it is abelian, and spanned over $\QQ_\ell$ by the classes constructed by A.Beilinson \cite{beilinson-modular} and in a different way by A.Beilinson and A.Levin \cite{beilinson-levin}. These classes are parametrized by a $\Ppp$-torsion point on $E$ and an element of the symmetric algebra of $H_1(E;\ZZ)$. These constructions work in the Hodge realization as well as in the $\ell$-adic one, and the mechanism of motivic correlators (described in \S2) gives alternative proofs of these statements. In particular, Beilinson and Levin's elliptic polylogarithms can be expressed in terms of the depth-1 Hodge correlator integrals -- Kronecker-Eisenstein series (\cite{beilinson-levin}, \S3). 

In this paper, we focus on the depth 2, the first case in which there is a nonzero Lie bracket. To describe the structure of the elliptic Galois Lie algebra, we can consider its standard cochain complex. Recall that the standard cochain complex of a Lie algebra $L$ is a complex of the exterior powers of its dual $L^\vee$, where the coboundary map $\delta$ is the dualization of the Lie bracket $[\,,\,]:L\wedge L\to L$:
\[{\rm CE}\ch(L^\vee)=\pq{0\to L^\vee\tto\delta L^\vee\wedge L^\vee\to L^\vee\wedge L^\vee\wedge L^\vee\to\dots}.\]
If $L^\vee$ is a graded Lie coalgebra, then ${\rm CE}\ch(L^\vee)$ is also graded. Applying the construction to the associated graded for the depth filtration of $\Lie_{(\ell)}^{\rm sym}(E,E[\Ppp])$, we obtain a cochain complex that is graded by weight and depth. In depth 1 the complex is concentrated in degree 1. However, the depth-2 part of this complex has a nontrivial coboundary map: the depth-2 elements map to wedge products of Beilinson-Levin classes:
\begin{equation}
    \gr^{D=2}\Lie^{\rm sym}_{(\ell)}(E,E[\Ppp])^\vee\to\pq{\gr^{D=1}\Lie^{\rm sym}_{(\ell)}(E,E[\Ppp])^\vee}^{\wedge2}.
    \label{eqn:d2_complex}
\end{equation}
We connect (the Hodge analogue of) this complex to the geometry of Bianchi hyperbolic threefolds.

\subsubsection*{Bianchi tesselation} Now let us describe the other side of the story. The Bianchi tesselation of the upper half-space $\HH^3$ for the ring $\ZZ[i]$ is the 3-dimensional version of the famous modular triangulation of the upper half-plane; the latter is the restriction of the Bianchi tesselation to the plane in $\HH^3$ lying above the real line (see Fig.~\ref{fig:bianchi}). This beautiful construction was given by L.Bianchi in 1892 \cite{bianchi}; see \cite{goncharov-euler} for a modern review. The fundamental domain is an octahedron with vertices at $0,1,i,i+1,\f{1+i}{2},\infty)$. Through the standard action of $\GL_2(\CC)$ on $\HH^3$, the group $\GL_2(\ZZ[i])$ acts transitively on the cells of the tesselation. If $\Ppp$ is a prime ideal in $\ZZ[i]$ and $\Gamma_1(\Ppp)\subset\GL_2(\ZZ[i])$ is the congruence subgroup \[\Gamma_1(\Ppp)=\cq{\begin{pmatrix}a&b\\c&d\end{pmatrix}\equiv\begin{pmatrix}1&0\\ * &1\end{pmatrix}\mod{\Ppp}},\]
the quotient $\Gamma_1(\Ppp)\sm\HH^3$ is a finite-volume hyperbolic manifold with cusps. 

\begin{figure}
    \includegraphics[height=3in]{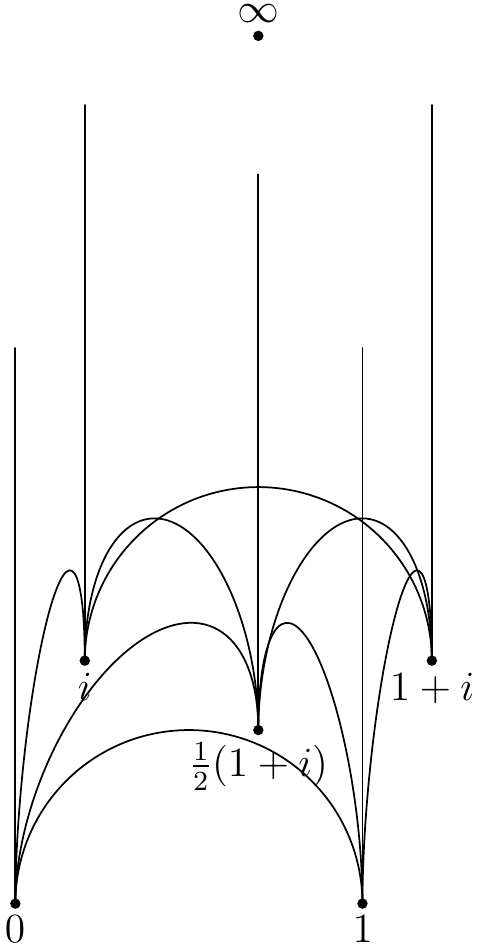}
    \hspace{1in}
    \includegraphics[height=3in]{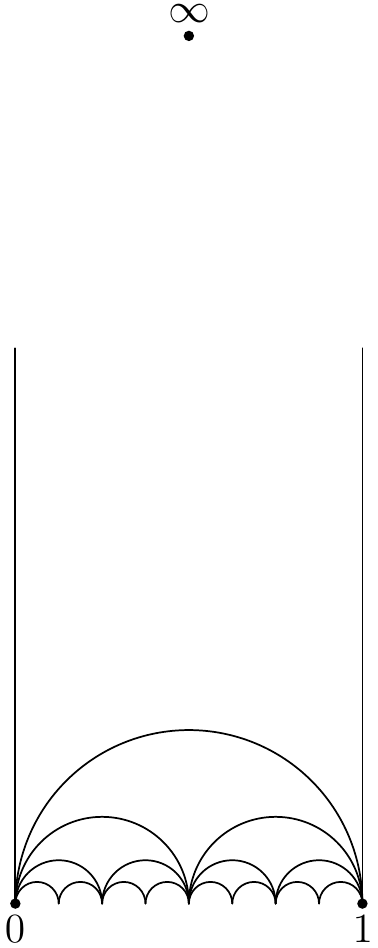}
    \simplecap{fig:bianchi}{\emph{Left:} The fundamental octahedron of the Bianchi tesselation for $\ZZ[i]$.\\\emph{Right:} The modular triangulation of the upper half-plane.}
\end{figure}

We build the following local system on this manifold. The group $H_1(E;\ZZ)$ has the structure of a $\ZZ[i]$-module, giving $H_1(E;\ZZ)\oplus H_1(E;\ZZ)$ the structure of a $\GL_2(\ZZ[i])$-module. We take its symmetric algebra $\Sym\ch(H_1(E;\ZZ)\oplus H_1(E;\ZZ))$. This $\GL_2(\ZZ[i])$-module determines an infinite-dimensional graded local system on $\Gamma_1(\Ppp)\sm\HH^3$. Denote this local system by $T_2$. 

Consider the chain complex of the Bianchi tesselation, placed in the cohomological degrees $[0,2]$. It is generated by the octahedral cells in degree 0, by ideal triangles in degree 1, and by geodesics in degree 2. Tensoring over $\Gamma_1(\Ppp)$ with $\Sym\ch(H_1(E;\ZZ)\oplus H_1(E;\ZZ))$, we get the chain complex of $\Gamma_1(\Ppp)\sm\HH^3$ with coefficients in the local system $T_2$.

\subsubsection*{The main construction} 

There is a Hodge analogue of $\Lie^{\rm sym}_{(\ell)}(E,E[\Ppp])$, denoted $\Lie^{\rm sym}_{\Hod}(E,E[\Ppp])$. In \S\ref{sec:relating} (Theorem~\ref{thm:bianchi_to_lie}), we construct a surjective morphism of complexes of graded $\ZZ[i]$-modules:
\begin{equation}
    \pq{\text{\parbox{0.45\textwidth}{\centering chain complex of the Bianchi orbifold $\Gamma_1(\Ppp)\sm\HH^3$ with coefficients in $T_2$}}}
\to\pq{\text{\parbox{0.3\textwidth}{\centering Hodge analogue of\\the complex (\ref{eqn:d2_complex})}}}
\label{eqn:main_morphism}
\end{equation}
In particular, we get surjective homomorphisms: \[H^i(\Gamma_1(\Ppp)\sm\HH^3,T_2)\to H^i(\gr^{D=2}\Lie_{\Hod}^{\rm sym}(E,E[\Ppp]),\QQ).\]
A key idea of A.Goncharov \cite{goncharov-euler}, which we develop further in this paper, is to map the cusps of the Bianchi orbifold $\Gamma_1(\Ppp)\sm\HH^3$ to $\Ppp$-torsion points of $E$.
This itself generalizes a similar picture for modular curves, first described in \cite{goncharov-manin}, where cusps of modular curves are identified with $p$-torsion points of $\GG_m$ in the study of the double logarithm at roots of unity.
When we advance to depth 2, the geodesics of the Bianchi tesselation map to wedge products of elements parametrized by $\Ppp$-torsion points, and the triangles must map to certain elements parametrized by three $\Ppp$-torsion points. 

Let us elaborate the map (\ref{eqn:main_morphism}) in each degree. In degree 2, we map a geodesic $(\alpha,\beta)$ of the Bianchi tesselation modulo $\Gamma_1(\Ppp)$, with a coefficient in the described local system to a wedge product of two depth-1 classes in the Galois Lie coalgebra. The data parametrizing a geodesic with a coefficient in $\Sym\ch(H_1(E;\ZZ)\oplus H_1(E;\ZZ))\cong\Sym\ch(H_1(E;\ZZ))^{\otimes2}$ is \emph{identical} to the data parametrizing a pair of classes in the image. In degree 1, the domain is generated by triangles in the Bianchi tesselation with a coefficient in the local system. Thus the image should be described in terms of elements depending on three $\Ppp$-torsion points and three elements of $\Sym\ch(H_1(E;\ZZ))$. These elements are \emph{motivic correlators}, which we introduce in the remainder of the introduction and in \S2. These elements can be visualized as a sequence of $\Ppp$-torsion points and 1-forms on $E$ written around a circle, modulo some relations. Their coproduct has a simple combinatorial description, and their real Hodge periods can be explicitly computed via Feynman integrals.

In summary, the maps are constructed as follows:
\begin{align*}
    \text{ideal triangle $(\alpha,\beta,\gamma)$}&\mapsto\text{element in $\gr^{D=2}\Lie^{\rm sym}(E,E[\Ppp])$ depending on $\Ppp$-torsion points $\alpha,\beta,\gamma$},\\
    \text{geodesic $(\alpha,\beta)$}&\mapsto\text{wedge product of Beilinson-Levin elements determined by $\alpha$ and $\beta$}.
\end{align*}
Remarkably, the combinatorial structure of the Bianchi tesselation is preserved in the space of motivic correlators, and thus the chain complex of the Bianchi orbifold \emph{maps surjectively} onto the standard cochain complex of a quotient of the Galois Lie algebra of $E\sm E[\Ppp]$.

\subsubsection*{Hodge realization} Now we will sketch this picture in the Hodge realization, which is the focus of this paper.

Let $E$ be a complex elliptic curve and $S\subset E$ a finite set of punctures. The pronilpotent completion of the fundamental group $\pi_1^\nil(E-S,v_0)$, with tangential base point at $v_0$, is a Lie algebra in the category of mixed $\QQ$-Hodge structures. The category of mixed $\QQ$-Hodge structures is canonically equivalent to the category of representations of a graded Lie algebra over $\QQ$. Let us take its image in the representation defining $\pi_1^\nil(E-S,v_0)$, and consider the graded dual Lie coalgebra $\Lie_\Hod^\vee(E,S)$.

The Hodge correlators, introduced by A.Goncharov in \cite{goncharov-hodge-correlators}, are canonical elements 
\begin{equation}
    \Cor_\Hod(\Omega_0,z_0,\dots,\Omega_n,z_n)\in\Lie_\Hod^\vee(E,S),
    \label{eqn:cor_hod}
\end{equation}
where $z_0,\dots,z_n\in S$ and $\Omega_1,\dots,\Omega_n$ are elements in the tensor algebra of $H^1(E;\CC)$. The coalgebra $\Lie_\Hod^\vee(E,S)$ carries a filtraion by depth; the element (\ref{eqn:cor_hod}) has depth $n$. These elements describe the real mixed Hodge structure on $\pi_1^\nil(E-S,v_0)\otimes\RR$. Their canonical real periods are the Hodge correlator functions, functions of $n+1$ points on $E$. We find new linear relations among the elements (\ref{eqn:cor_hod}).

At a cusp on the modular curve, as $E$ degenerates to the nodal projective line, these relations specialize to known relations among periods of the mixed Tate motive associated with $\PP^1$ punctured at a finite set of points. If $n=2$, our elliptic relations specialize to the full set of \emph{double shuffle relations}, the most general known relations, which were previously described by the author using Hodge correlators (\cite{malkin-shuffle}).

Suppose that $E$ is one of the CM elliptic curves $\CC/(\ZZ+\ZZ i)$ or $\CC/(\ZZ+\ZZ\pq{\f{1+\sqrt{-3}}{2}})$, $\OOO=\End E$, and $\Ppp$ is a prime in $\OOO$. The subalgebra $\Lie_\Hod^{\rm sym}(E,E[\Ppp])$ of $\Lie_\Hod(E,S)$ is constructed as in the $\ell$-adic case. We construct the morphism (\ref{eqn:main_morphism}) in this setting, where the object standing on the right is the complex ${\rm CE}\ch\pq{\gr^{D=2}\Lie_\Hod^{\rm sym}(E,E[\Ppp])}$.

Our construction simultaneously generalizes several results of A.Goncharov:
\begin{enumerate}
\item \emph{The relation between Voronoi complexes and mixed Tate motives: } The Bianchi complexes are the higher-degree analogues of the Voronoi complexes, complexes of $\GL_k(\ZZ)$-modules from tesselations of the upper half-plane $\HH^2$. A map from the Voronoi complexes to motivic objects associated with rational curves punctured at roots of unity constructed for $k=2,3,4$, using either multiple polylogarithms (\cite{goncharov-polylogs-modular}) or motivic correlators (\cite{goncharov-motivic-modular}), which satisfy the double shuffle relations. The relations we found for elliptic motivic correlators in depth 2 are deformations of the second shuffle relations.
\item \emph{Euler complexes: } The map from the Bianchi complexes to a space of motivic theta functions on elliptic curves constructed by \cite{goncharov-euler} in depth 2 and weight 4. We generalize this construction to all weights: \cite{goncharov-euler}'s map is the restriction of our map to the trivial local system.
\end{enumerate}

\subsubsection*{Structure}

In \S\ref{sec:hodge_and_motivic}, we review the construction of Hodge correlators. We in particular explain our results on the level of Hodge correlator integrals.

In \S\ref{sec:mc_elliptic_curves} we establish some properties of motivic correlators on elliptic curves. The main new result of this section is the dihedral symmetry relation for depth 2 correlators (Theorem~\ref{thm:dihedral_depth2}).

In \S\ref{sec:bianchi_and_modular} we review the definitions of the Bianchi complexes, define the modular complexes for imaginary quadratic fields, and construct a map between the two in the Gaussian and Eisenstein cases. In \S\ref{sec:relating} we combine the results of the two preceding sections to prove the main results relating Bianchi complexes and the elliptic Galois Lie algebra.

In \S\ref{sec:application}, we show how our results generalize those of \cite{goncharov-polylogs-modular,goncharov-euler,malkin-shuffle}.

\subsubsection*{Acknowledgements}

The author is grateful to A.B.\ Goncharov for suggesting this problem, for many helpful explanations, and for comments on a draft of this paper.

This material is based in part upon work supported by NSF grants 
DMS-1440140, 1107452, 1107263, and 1107367.

\section{Hodge and motivic correlators}

\label{sec:hodge_and_motivic}

\subsection{Real Hodge point of view: Relations on Hodge correlator integrals}

Let $X$ be a complex curve (in this paper, $X=\PP^1(\CC)$ or an elliptic curve). The Hodge correlator functions, defined in \cite{goncharov-hodge-correlators}, are functions
\[\Cor_\HHH(x_0,\dots,x_n),\] 
where each $x_i$ is either a point of $X$ or a 1-form representing a class in $H^1(X;\CC)$. The \emph{depth} of this expression is the number of points among the $x_i$ minus one, if $X$ is an elliptic curve, or the number of nonzero points among the $x_i$ minus one, if $X=\PP^1$. The \emph{weight} is $n$ plus the depth.

The Hodge correlators depend on a choice of a base point $s\in X$ and a tangent vector $v_0$ at $s$. If $n=1$ and $x_0,x_1\in X$, then $\Cor_\HHH(x_0,x_1)$ is a (normalized) Green's function with pole at $s$. In particular,
\begin{itemize}
\item If $X=\PP^1$ and $s=\infty$, then \[\Cor_\HHH(x_0,x_1)=G_\infty(x_0,x_1)=(2\pi i)\inv\log\aq{x_0-x_1}+C.\] The constant $C$ depends on the choice of tangent vector at $\infty$, but the correlator is independent of this constant in weight $>2$, so we will ignore it when convenient. The correlator for other tangential base points can be derived using the fact that it is invariant under automorphisms of $\PP^1$ acting on the base point and the arguments.
\item If $X$ is the elliptic curve $\CC/(\ZZ+\ZZ\tau)$, where $\Im(\tau)>0$, with coordinate $z$ inherited from the complex plane, then \[\Cor_\HHH(x_0,x_1)=G_s(x_0,x_1)=G_{\rm Ar}(x_0,x_1)-G_{\rm Ar}(x_0,s)-G_{\rm Ar}(s,x_1)+C.\] Here $G_{\rm Ar}$ is the Arakelov Green's function, the unique solution to the elliptic partial differential equation $(2\pi i)\inv\del\ol\del G_{\rm Ar}(x)={\rm vol}_E-\delta_0$. It has the Fourier expansion 
\begin{equation}
    G_{\rm Ar}(z)=\f{2\Im(\tau)}{2\pi i}\sum_{\gamma\in(\ZZ+\ZZ\tau)\sm\cq0}\f{\exp\pq{2\pi i\Im(z\ol\gamma)/\Im(\tau)}}{\aq\gamma^2}.
    \label{eqn:fourier}
\end{equation}
The Arakelov Green's function has a logarithmic singularity at 0. Hence, the function $\Cor_\HHH(x_0,x_1)$ has singularities of the form $\log\aq z$ at the divisors $x_0=x_1$, $x_0=s$, $x_1=s$.
\item \textbf{Remark: } The Green's function on $\PP^1$ is a specialization of the one on $E$. Precisely, write $G^{E_\tau}$ for the Green's function on $E=\CC/(\ZZ+\ZZ\tau)$ with base point 0. Then, taking $z$ to be the coordinate on $E_\tau$ inherited from the complex plane, such that the section $z\in E_\tau$ approaches $e^{2\pi iz}\in\PP^1$, with appropriate choice of tangential base points, \[\lim_{\tau\to+i\infty}G^{E_\tau}(z_1,z_2)=G_{1}\pq{e^{2\pi iz_1},e^{2\pi iz_2}}=\log\aq{\f{1}{e^{2\pi iz_1}-1}-\f{1}{e^{2\pi iz_2}-1}}.\] (This can be shown by a residue computation or an application of the Kronecker limit formula. We will require this fact in \S\ref{sec:degen}.)
\end{itemize}

If $n\geq2$, the Hodge correlators are defined as a sum of integrals depending on plane trivalent trees. Picture the $x_0,\dots,x_n$ written counterclockwise along the boundary of a disc, and consider a trivalent tree $T$ embedded in the disc with leaves at the $n+1$ boundary points. The tree has $n-1$ interior vertices $V^\circ$ and $2n-1$ edges $E_0,\dots,E_{2n-2}$. The embedding into the plane gives a canonical orientation ${\rm Or}_T\in\cq{\pm1}$ (an ordering of the edges up to even permutation).

Assign to each interior vertex $v\in V^\circ$ a copy of $X$, called $X_v$, with coordinate $x_v$. Then assign to each edge $E_i$ either a function $f_i$ or a 1-form $\omega_i$, as follows:
\begin{enumerate}[(1)]
\item If $E_i=(u,v)$ is an interior edge, let $f_i=G_s(x_u,x_v)$, a function on $X_u\times X_v$.
\item If $E_i=(u,x_j)$ is a boundary edge with the leaf decorated by a point $x_j\in X$, let $f_i=G_s(x_u,x_j)$, a function on $X_u$.
\item If $E_i=(u,x_j)$ is a boundary edge with $x_j=\omega$ a 1-form, let $\omega_i=\omega(x_u)$, a 1-form on $X_u$.
\end{enumerate}
Without loss of generality, $E_0,\dots,E_k$ are the edges labeled by a function (i.e., not boundary edges decorated by a 1-form). Suppose also that each form is either purely holomorphic or purely antiholomorphic (which we may do because the Hodge correlators are linear in the forms); let there be $p$ and $q$ forms of these types, respectively. Then, setting $d^\CC=\del-\ol\del$, we define
\begin{equation}
c_T(x_0,\dots,x_n)=(-2)^k\binom{k}{\f12(k+p-q)}\inv{\rm Or}_T\int_{X^{V^\circ}}f_0\,d^\CC f_1\wedge\dots\wedge d^\CC f_k\wedge\omega_{k+1}\wedge\dots\wedge\omega_{2n-2}.
\end{equation}
The Hodge correlator is the sum of such expressions over all plane trivalent trees,
\begin{equation}
\Cor_\HHH(x_0,\dots,x_n)=\sum_Tc_T(x_0,\dots,x_n).\label{eqn:def_hodge_corr}
\end{equation}
The Hodge correlator is independent of the choice of ordering of edges. As a function of the arguments that are points on $X$, it is either purely real or purely imaginary.

Fig.~\ref{fig:hc_example} shows a simple example of the integral corresponding to one of the two trees contributing to $\Cor_\HHH(a,b,c,\o)$.

\begin{figure}
    \centering
    \includegraphics[width=0.2\textwidth]{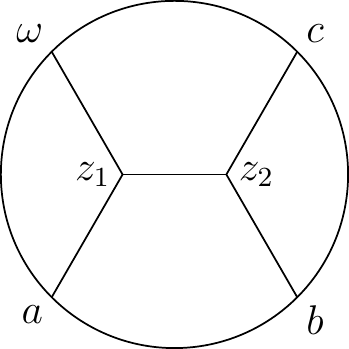}

    \[\int_{z_1,z_2}G(z_1,a)\,d^\CC G(z_1,z_2)\wedge d^\CC G(z_2,b)\wedge d^\CC G(z_2,c)\wedge\omega(z_1)\]
\simplecap{fig:hc_example}{One of the trees contributing to $\Cor_\HHH(a,b,c,\o)$.}
\end{figure}

The Hodge correlators satisfy a family of \emph{(first) shuffle relations}. For $i,j>0$, let $\Sigma_{i,j}$ be the set of $(i,j)$-shuffles, permutations $\sigma\in S_{i+j}$ such that $\sigma(1)<\dots<\sigma(i)$ and $\sigma(i+1)<\dots<\sigma(i+j)$. The $(i,j)$-shuffle relation states:
\begin{equation}
\sum_{\sigma\in\Sigma_{i,j}}\Cor_\HHH(x_0,x_{\sigma\inv(1)},x_{\sigma\inv(2)},\dots,x_{\sigma\inv(i+j)})=0.
\end{equation}
For Hodge correlators of depth 2 on an elilptic curve with arbitrary base point, we found a second shuffle relation. It has the form:
\begin{equation}
\Cor_\HHH(S_{n_0,n_0'},0,S_{n_1,n_1'},a,S_{n_2,n_2'},a+b)
+\Cor_\HHH(S_{n_0,n_0'},0,S_{n_2,n_2'},b,S_{n_1,n_1'},a+b)
+\text{lower-depth terms}=0,
\label{eqn:ec_second_shuffle}
\end{equation}
where an argument $S_{n,n'}$ indicates that we sum over all possible ways to insert in some order the arguments \[\underbrace{\o,\dots,\o}_n,\underbrace{\ob,\dots,\ob}_{n'}.\] 
The highest-depth terms in these relations arise from shuffles of the \emph{differences} between successive arguments, $x_i-x_{i-1}$, together with the 1-forms between those arguments. For example, in (\ref{eqn:ec_second_shuffle}) we have shuffled $a$ (with $n_1$ copies of $\o$ and $n_1'$ of $\ob$) with $b$ (with $n_2$ $\o$'s and $n_2'$ $\ob$'s).

We describe the lower-depth correction terms in \S\ref{sec:mc_modulo_depth}. In the simplest case -- weight 4 -- the full relation is:
\begin{align*}
\Cor_\HHH(0,a,a+b)
+\Cor_\HHH(0,b,a+b)
-(&\Cor_\HHH(0,\o,\ob,a+b)
+\Cor_\HHH(0,\ob,\o,a+b))\\
-\f12\biggl(
&\Cor_\HHH(0,a,\o,\ob)
-\Cor_\HHH(0,a,\ob,\o)\\
&+\Cor_\HHH(0,b,\o,\ob)
-\Cor_\HHH(0,b,\ob,\o)\\
&+\Cor_\HHH(\o,\ob,a,a+b)
-\Cor_\HHH(\ob,\o,a,a+b)\\
&+\Cor_\HHH(\o,\ob,b,a+b)
-\Cor_\HHH(\ob,\o,b,a+b)
\biggr)&=0.
\end{align*}
The relation (\ref{eqn:ec_second_shuffle}) can be formulated simply as a functional equation on biperiodic functions of three complex variables. It states that a sum of several integrals over an elliptic curve is equal to 0, modulo correlators of depth 1 (which are expressed by Kronecker-Eisenstein series). However, this functional equation is difficult to prove. To understand it, we will need to upgrade it to the Hodge-theoretic or motivic setting.

\begin{figure}
    \centering
    \includegraphics[height=0.2\textwidth]{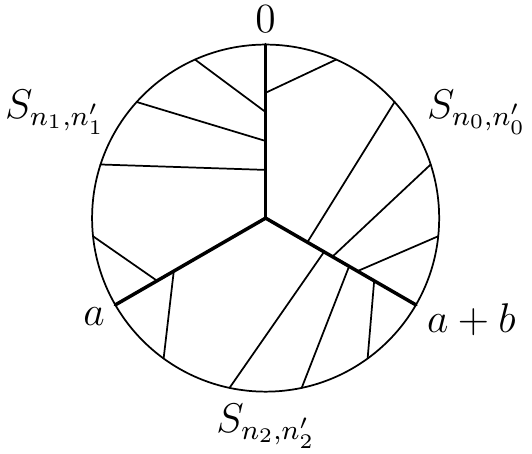}
    \parbox{0.1\textwidth}{\centering\vskip-0.1\textwidth\Huge+\vskip0.1\textwidth}
    \includegraphics[height=0.2\textwidth]{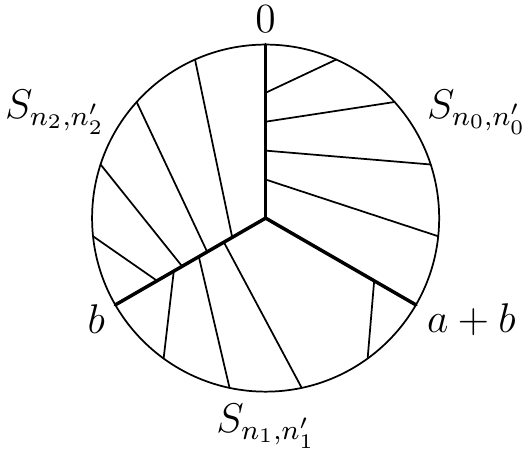}
\simplecap{fig:second_shuffle}{The highest-depth terms of the second shuffle relation on an elliptic curve.}
\end{figure}

The second shuffle relations have a prehistory. The first objects known to satisfy a system first and second shuffle relations of this form were the multiple polylogarithms (see \cite{goncharov-polylogs-modular}). These relations follow from two alternative expressions for multiple polylogarithms: as power series and as iterated integrals. In \cite{malkin-shuffle}, for $X=\PP^1$, the author found second shuffle relations for Hodge correlators, in every depth, and described the lower-depth terms. In depth 2, these relations depend on integers $n_0,n_1,n_2\geq0$ and points $a,b\in\GG_m\sm\cq1$. They state:
\begin{align}
\Cor_\HHH(\underbrace{0,\dots,0}_{n_0},1,\underbrace{0,\dots,0}_{n_1},a,\underbrace{0,\dots,0}_{n_2},ab)
+\Cor_\HHH(\underbrace{0,\dots,0}_{n_0},1,\underbrace{0,\dots,0}_{n_2},b,\underbrace{0,\dots,0}_{n_1},ab)&\nonumber\\
+\,\text{lower-depth terms}&=0.
\label{eqn:p1_second_shuffle}
\end{align}
The highest-depth terms in these relations arise from shuffles of the quotients between successive arguments, $\f{x_i}{x_{i-1}}$, together with the 0s between those arguments. For example, in (\ref{eqn:p1_second_shuffle}) we have shuffled $a$ (with $n_1$ 0s) with $b$ (with $n_2$ 0s).  See Fig.~\ref{fig:second_shuffle}) for an illustration.

Conjecturally, the first and second shuffle relations give all linear relations among the Hodge correlators on $\PP^1$. While the first shuffle relations emerge from the trivalent tree construction -- they hold on the level of the \emph{integrands} in (\ref{eqn:def_hodge_corr}) -- the proof of the second shuffle relations is difficult, requiring motivic or Hodge-theoretic arguments even in depth 2. %

Note the similarity between (\ref{eqn:p1_second_shuffle}) and (\ref{eqn:ec_second_shuffle}). In fact, as an elliptic curve degenerates to a nodal projective line, a variant of the second shuffle relation (\ref{eqn:ec_second_shuffle}) specializes to (\ref{eqn:p1_second_shuffle}).

\subsection{Hodge-theoretic / motivic point of view: Correlators and motivic $\pi_1$}

We briefly review the construction of Hodge and motivic correlators from \cite{goncharov-hodge-correlators}. Hodge correlators are objects in the fundamental Lie coalgebra of the category of $\RR$-mixed Hodge structures, and are Hodge-theoretic upgrades of the Hodge correlator functions.

\subsubsection{Summary}

In \cite{goncharov-hodge-correlators}, the Hodge correlator functions $\Cor_\HHH(x_0,\dots,x_n)$ of the previous section were upgraded to elements of the Tannakian Lie coalgebra $\Lie_\Hod^\vee$ of the category of real mixed Hodge structures
\begin{equation}
    \Cor_\Hod(z_0,\dots,z_n)\in\Lie_\Hod^\vee.\label{eqn:elem_cor_hod}
\end{equation} 
The upgraded Hodge correlators (\ref{eqn:elem_cor_hod}) that satisfy the first shuffle relations, and their coproduct in the coalgebra $\Lie_\Hod^\vee$ is given by a simple formula, which we give below.

One of the main results of this paper is that the elements (\ref{eqn:elem_cor_hod}) satisfy  a second shuffle relation in depth 2.

\subsubsection{Hodge-theoretic setup}

Let $\MHS_\RR$ of be the tensor category of $\RR$-mixed Hodge structures and $\HS_\RR$ the category of $\RR$-pure Hodge structures. Every object of $\MHS_\RR$ is filtered by weight, and  $\MHS_\RR$ is generated by the simple objects $\RR(p,q)+\RR(q,p)$ ($p,q\in\ZZ$). By Deligne's theory \cite{deligne-hodge3}, the cohomology of a (possibly singular) complex variety is a mixed Hodge structure.

The \emph{Galois Lie algebra} of the category of mixed Hodge structures, $\Lie_\Hod$, is the algebra of tensor derivations of the functor $\gr^W:\MHS_\RR\to\HS_\RR$. It is a graded Lie algebra in the category $\HS_\RR$, and $\MHS_\RR$ is equivalent to the category of graded $\Lie_\Hod$-modules in $\HS_\RR$. Let $\Lie_\Hod^\vee$ be its graded dual. A canonical \emph{period map} \[p:\Lie_\Hod^\vee\to\RR\] was defined in \cite{goncharov-hodge-correlators}, \S1.11.

Let $X$ be a smooth curve, $S\subset X$ a finite set of punctures, $s\in S$ a distinguished puncture (called the base point), and $v_0$ a distinguished tangent vector at $s$. The pronilpotent completion $\pi_1^\nil(X\sm(S\cup\cq s),v_0)$ of the fundamental group $\pi_1(X\sm S,s)$ carries a mixed Hodge structure, depending on $v_0$, and thus there is a map
\[\Lie_\Hod\to\Der\pq{\gr^W\pi_1^\nil(X\sm S,v_0)}.\]

\subsubsection{Hodge correlator coalgebra}

\label{sec:hodge_cor_coalg}

The \emph{Hodge correlator coalgebra} is defined by \cite{goncharov-hodge-correlators} as 
\begin{equation*}
\CLie_{X,S,v_0}^\vee:=\f{T(\CC\bq{S}^\vee\oplus H^1(X;\CC))}{\text{relations}}\otimes H_2(X).
\end{equation*}
Note that $H_2(X)\cong\RR(1)$. If $[h]\in H_2(X)$ is the fundamental class, we write $x(1)$ for $x\otimes[h]$. This coalgebra is graded by weight. It is more finely graded by the Hodge bidegree, or \emph{type}, where points in $S$ have type $(1,1)$, holomorphic and antiholomorphic 1-forms have type $(1,0)$ or $(0,1)$, respectively, and $H_2(X)$ has type $(-1,-1)$, extended to be additive with respect to the tensor product. The weight of an element of type $(p,q)$ is $p+q$.

The relations are the following:
\begin{enumerate}[(1)]
\item Cyclic symmetry: $x_0\otimes\dots\otimes x_n=x_1\otimes\dots\otimes x_n\otimes x_0$.
\item (First) shuffle relations:
\begin{equation*}
\sum_{\sigma\in\Sigma_{i,j}}x_0\otimes x_{\sigma\inv(1)}\otimes\dots\otimes x_{\sigma\inv(i+j)}=0.
\end{equation*}
\item Take the quotient by the elements of nonpositive weight.
\end{enumerate}
An action of the graded dual Lie algebra $\CLie_{X,S,v_0}$ by derivations on $\gr^W\pi_1^\nil(X\sm S,v_0)\otimes\CC$
was constructed by \cite{goncharov-hodge-correlators}. This action is injective; its image consists of the \emph{special derivations} \[\Der^S\pq{\gr^W\pi_1^\nil(X\sm S,v_0)\otimes\CC},\] those which act by 0 on the loop around $\infty$ and preserve the conjugacy classes of all the loops $s\in S\sm\cq{s}$.

Dualizing this map composed with the action of $\Lie_\Hod$, we get the \emph{Hodge correlator morphism} of Lie coalgebras:
\[\Cor_\Hod:\CLie^\vee_{X,S,v_0}\to\Lie_\Hod^\vee.\]
Let $\Lie_\Hod^\vee(X,S,v_0)$ denote the image of this action, and let $\Lie_\Hod^\vee(X,S)$ denote the algebra generated by the $\Lie_\Hod^\vee(X,S,v_0)$ for all choices of base point. (Below, we will fix $S=E[p]$ for $E$ an elliptic curve, so $\Lie_\Hod^\vee(X,S)$ does not depend on the choice of base point in $S$.)
We will also write $\Cor_\Hod(x_0,\dots,x_n)$ for $\Cor_\Hod\pq{(x_0\otimes\dots\otimes x_n)(1)}$, or $\Cor_s(\dots)$, when we wish to specify the base point.

The Lie coalgebra structure on $\CLie_{X,S,v_0}^\vee$ has a simple description on the generators. There are two terms in the coproduct, $\delta_{\rm D}$ and $\delta_{\rm Cas}$, which are each sums over ``cuts'' of the element \[C=\pq{x_0\otimes\dots\otimes x_n}\otimes[h],\] which we picture as $x_0,\dots,x_n$ written counterclockwise around a circle. 

\begin{enumerate}[(1)]
\item Term $\delta_S$: Consider a line inside the circle beginning at a point on the circle labeled by a puncture $x_i$ and ending between two adjacent points. It cuts the circle into two parts $C_1$ and $C_2$, which share only the point $x_i$, where $C_1$ lies clockwise of $x_i$. This contributes to the coproduct the term $C_1\wedge C_2$, and $\delta_{\rm D}C$ is the sum of these terms over all such cuts. That is,
\[\delta_SC=\sum_{\stackrel{\rm cyc}{x_0\in S}}\sum_{p=1}^n\pq{\pq{x_0\otimes x_p\otimes\dots\otimes x_n}\otimes[h]}\wedge\pq{\pq{x_0\otimes x_1\otimes\dots\otimes x_{p-1}}\otimes[h]},\] where the outer sum is only taken over those cyclic reorderings where $v_0$ is a puncture.
(See Fig.~\ref{fig:coproduct}, top.)
\item Term $\delta_{\rm Cas}$: Consider a line inside the circle beginning between two points $y_1$ and $z_1$ and ending between two points $y_2$ and $z_2$. It cuts the circle into two parts $C_1$ and $C_2$, in which $y_1$ and $z_2$ are adjacent and in which $y_2$ and $z_1$ are adjacent. We insert a point labeled $\omega$ between $y_1$ and $z_2$ on $C_1$ and a point labeled $\omega^\vee$ between $z_2$ and $y_1$ on $C_2$ to obtain $C_1'$ and $C_2'$, then take the sum over $\omega$ in a fixed symplectic basis $\cq{\omega_i}$ of $H^1(X;\CC)$. This contributes the term $C_1'\wedge C_2'$, and $\delta_{\rm Cas}$ is the sum of these terms over all such cuts. That is,
\begin{equation}
    \delta_{\rm Cas}C=\sum_{p=0}^n\sum_{q=0}^n\sum_{i=1,2}\pq{\pq{x_p\otimes\dots\otimes x_{q-1}\otimes\omega_i}\otimes[h]}\wedge\otimes\pq{\pq{x_q\otimes\dots\otimes x_{p-1}\otimes\omega_i^\vee}\otimes[h]}.
\end{equation}
(See Fig.~\ref{fig:coproduct}, bottom.)
\end{enumerate}
The term $\delta_{\rm Cas}$ are absent if $X=\PP^1$. If $E$ is an elliptic curve, $\CLie_{X,S,v_0}^\vee$ is graded by weight and filtered by depth, and the terms $\delta_{\rm Cas}$ disappear in the associated graded $\gr^D\CLie_{X,S,v_0}^\vee$.

\begin{figure}[ht]
\centering
\includegraphics[width=0.7\textwidth]{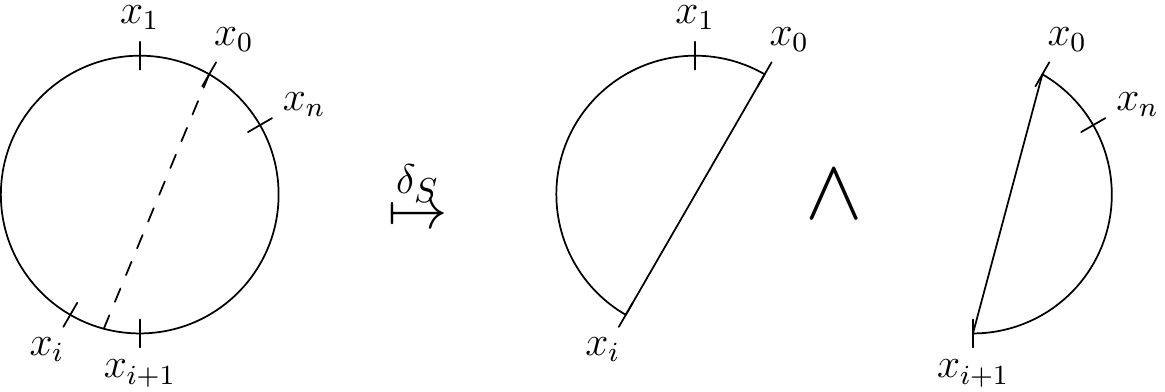}\\
\includegraphics[width=0.7\textwidth]{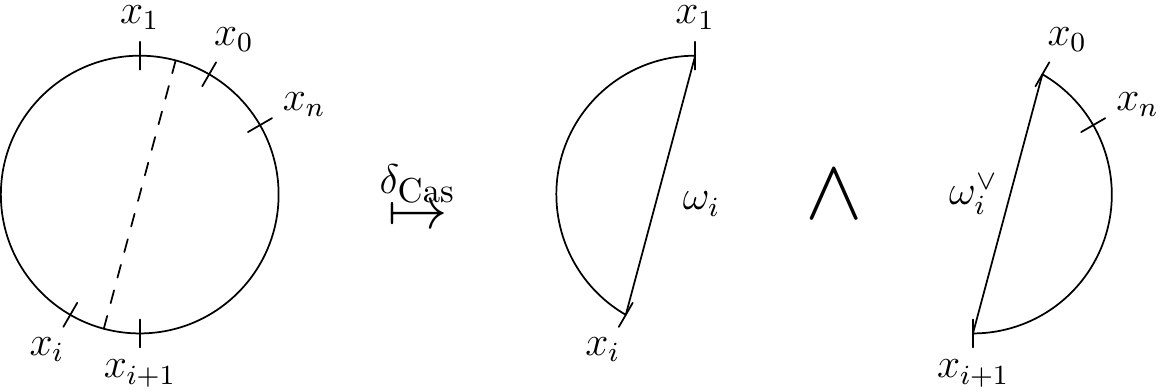}
\simplecap{fig:coproduct}{\emph{Above:} The typical term in the $\delta_S$ component of the coproduct.\\\emph{Below:} The typical term in the $\delta_{\rm Cas}$ component.}
\end{figure}

\subsubsection{Period map and Hodge correlator functions}

Recall that the Hodge correlator functions $\Cor_\HHH(x_0,\dots,x_n)$ satisfy cyclic symmetry and shuffle relations, so we may also denote by $\Cor_\HHH$ the function
\begin{align*}
    \Cor_\HHH:\CLie^\vee_{X,S,v_0}&\to\CC,\\
    (x_0\otimes\dots\otimes x_n)(1)&\mapsto\Cor_\HHH(x_0,\dots,x_n).
\end{align*}
The dual to the Hodge correlator $\Cor_\HHH:\CLie_{X,S,v_0}^\vee\to\CC$, an element of $\CLie_{X,S,v_0}$, is called the \emph{Green operator} $\mathbf{G}_{v_0}$. It can be viewed as a special derivation of $\gr^W\pi_1^\nil(X\sm S,v_0)\otimes\CC$, and defines a real mixed Hodge structure on $\pi_1^\nil(X\sm S,v_0)$. An element $x\in\CLie^\vee_{X,S,v_0}$ of type $(p,q)$ provides a framing $\RR(p,q)+\RR(q,p)\to\gr^W_{p+q}\pi_1^\nil(X\sm S,v_0)$, and $\Cor_\Hod(x)$ is the element of $\Lie_\Hod^\vee$ induced by this framing.

As made precise by a main result of \cite{goncharov-hodge-correlators}, $\Cor_\HHH$ factors through the Hodge correlator map to $\Lie_\Hod^\vee$ and the period map $\Lie_\Hod^\vee\to\CC$, and the resulting mixed Hodge structure on $\pi_1^\nil$ coincides with the standard one.

\begin{thm}[\cite{goncharov-hodge-correlators}, Theorem 1.12]
    \begin{enumerate}[(a)]
    \item For $x\in\CLie^\vee_{X,S,v_0}$, $\Cor_\HHH(x)=(2\pi i)^{-n}p(\Cor_\Hod(x))$, where $p$ is the canonical period map $\Lie_\Hod^\vee\to\RR$.
    \item The mixed Hodge structure on $\pi_1^\nil$ determined by the dual Hodge correlator map coincides with the standard mixed Hodge structure on $\pi_1^\nil$.
    \end{enumerate}
    \label{thm:hc_main_point}
\end{thm}

Furthermore, let $X/B$ be a smooth curve over a base $B$. For a collection of nonintersecting sections $S$ and choice of relative tangent vector $v_0$, we can analogously define $\CLie^\vee_{X/B,S,v_0}$. In this setting, for $x\in\CLie^\vee_{X/B,S,v_0}$, \cite{goncharov-hodge-correlators} constructs a connection on the fiberwise $\Cor_\Hod(x)$ that makes $\Cor_\Hod(x)$ a variation of mixed Hodge structures over $B$. We have the following essential fact, which follows from the Griffiths transversality condition:
\begin{lma}
If $x\in\CLie^\vee_{X/B,S,v_0}$ of type $(p,q)$, and weight $p+q=n>2$, has $\delta(\Cor_\Hod(x))=0$, and $\Cor_\HHH(x|_b)=0$ at some $b\in B$, then $\Cor_\Hod(x)=0$.
\label{lma:d0rigid}
\end{lma}
\begin{proof}
If $\delta(\Cor_\Hod(x))=0$, then $\Cor_\Hod(x)\in\Ext^1(\RR(0),\RR(p,q)+\RR(q,p))$, which is one-dimensional and rigid by the Griffiths transversality condition. Hence the variation is constant and captured by the period $p:\Ext^1(\RR(0),\RR(p,q)+\RR(q,p))\to\CC$.
\end{proof}

One of the main results of this paper is that the relations (\ref{eqn:ec_second_shuffle}) hold for the elements $\Cor_\Hod$: the equality between functions is upgraded to a relation in the fundamental Lie coalgebra of mixed Hodge structures.

\subsubsection{Motivic correlators}

Let $F$ be a number field. Beilinson's conjectures (\cite{beilinson-height-pairing}) predict that there is a category $\MM_F$ of mixed motives over $F$. Every object in $\MM_F$ should have a weight filtration, and there should be a functor $\gr^W:\MM_F\to\PM_F$, where $\PM_F$ is the category of pure motives over $F$. For every embedding $\sigma:F\to\CC$, there should be a realization functor $r_\sigma:\MM_F\to\MHS_F$. For every simple object $M\in\MM_F$, there should be an injective regulator map
\[\reg:\Ext_{\MM_F}^1(\QQ(0),M)\to\bigoplus_{F\to\CC/\text{complex conj.}}\Ext^1_{\MHS_\RR}(\RR(0),r_\sigma(M))).\]
The \emph{fundamental (motivic) Lie algebra} $\Lie_{\Mot/F}$ is the algebra of tensor derivations of the functor $\gr^W$, a graded Lie algebra in the category $\PM_F$, and $\MM_F$ is equivalent to the category of graded $\Lie_{\Mot/F}$-modules. An embedding $\sigma$ induces a map $r_\sigma:\Lie_{\Mot/F}^\vee\to\Lie_\Hod^\vee$.

Let $X$ be a curve defined over $F$, $S\subset X(F)$ a finite set of punctures, and $v_0$ the distinguished tangent vector at $s\in S$. There is expected to be a \emph{motivic fundamental group} $\pi_1^\Mot(X\sm S,v_\infty)$, a prounipotent group scheme in the category $\MM_F$. The Hodge realization of its Lie algebra should be $\pi_1^\nil(X\sm S,v_0)$. As it is an object in $\MM_F$, there is an action $\Lie_{\Mot/F}\to\Der\pq{\gr^W\pi_1^\Mot}$.

The construction of the Hodge correlator coalgebra $\CLie^\vee_{X,S,v_0}$ can be upgraded to the motivic setting, simply by replacing all the Hodge-theoretic objects by their motivic avatars.
For example, the definition of the \emph{motivic correlator coalgebra} mimics that of its Hodge realization:
\begin{equation*}
\pq{\CLie^\Mot_{X,S,v_0}}^\vee:=\f{T\pq{(\QQ(1)^{S})^\vee\oplus H^1(X)}}{\text{relations}}\otimes H_2(X),
\end{equation*}
a graded Lie coalgebra in the category of pure motives over $F$, where the relations imposed are the cyclic symmetry, first shuffles, and quotient by nonpositive weight. Then $\CLie^\Mot_{X,S,v_0}$ is isomorphic to the algebra of special derivations of $\gr^W\pi_1^\Mot(X-S,v_0)$, and there is a \emph{motivic correlator map}
\begin{equation*}
\Cor_\Mot:\pq{\CLie^\Mot_{X,S,v_0}}^\vee\to\Lie_{\Mot/F}^\vee.
\end{equation*}
We will write $\Cor_\Mot(x_0,\dots,x_n)$ for $\Cor_\Mot(\pq{x_0\otimes\dots\otimes x_n}(1))$.

Fix an embedding $\sigma:F\to\CC$. We have the composition of the realization functor with the period map:
\[
    \Cor_\HHH\circ r_\sigma:\pq{\CLie^\Mot_{X,S,v_0}}^\vee\otimes\CC\to\CLie_{X,S,v_0}^\vee\otimes\CC\to\CC.
\]

By Theorem~\ref{thm:hc_main_point}, it coincides with the composition 
\[
    \pq{\CLie^\Mot_{X,S,v_\infty}}^\vee\to\Lie_\Mot^\vee\to\Lie_\Hod^\vee\to\CC.   
\]
We can summarize all of the described objects and maps defined as follows:
\[
\xymatrix{
    \Der^S(\gr^W\pi_1^\Mot(X\sm S,v_0))^\vee\ar@{-}[r]
    &(\CLie_{X,S,v_0}^\Mot)^\vee\ar[r]^{\quad\Cor_\Mot}\ar[d]^r
    &\Lie_{\Mot/F}^\vee\ar[d]^r
    \\
    \Der^S(\gr^W\pi_1^\nil(X\sm S,v_0))^\vee\ar@{-}[r]
    &(\CLie_{X,S,v_0}^\vee)\ar[r]^{\quad\Cor_\Hod}\ar[dr]_{\Cor_\HHH}&\Lie_\Hod^\vee\ar[d]^p
    \\
    &&\CC.
}    
\]

Relations among the motivic correlators can be proven by showing that they hold in the Hodge realization under any complex embedding. Precisely, there is the following fact, which is an immediate consequence of the (hypothetical) injectivity of the regulator and Lemma~\ref{lma:d0rigid}.
\begin{lma}
    Suppose $x\in\pq{\CLie_{X,S,v_\infty}^\Mot}^\vee$ is of type $(p,q)$ with weight $p+q>2$, $\delta\Cor_\Mot(x)=0$, and $\Cor_\HHH(r(x))=0$ for every embedding $r:F\to\CC$. Then $\Cor_\Mot(x)=0$.
    \label{lma:d0h0_rational}
\end{lma}
This fact allows us to lift relations on Hodge correlators to relations on motivic correlators. In particular, all results in this paper -- the second shuffle relations for Hodge correlators and the map from the Bianchi complexes to an algebra of Hodge correlators -- should hold with ``Hodge'' replaced by ``motivic''.

Assuming the motivic formalism, the results in the Hodge realization can then be translated to the $\ell$-adic realization, via the motivic correlators. In particular, the results stated in the introduction would hold for the $\ell$-adic elliptic Galois algebra.

\section{Motivic correlators on elliptic curves}

\label{sec:mc_elliptic_curves}

\subsection{Main properties}

\subsubsection{Definitions}

We work with a complex elliptic curve $E$. Recall $S\subset E$ is a finite set of punctures. Let $\OOO=\End(E)$, so either $\OOO=\ZZ$ or a lattice in an imaginary quadratic field.

Let $\omega,\ol\omega$ be a symplectic basis for $H^1(E;\CC)$. $\CL_{E,S,v_0}^\vee$ is generated by elements
\begin{align*}
    C_s(\Omega_0,s_0,\dots,\Omega_n,s_n)
    =&\underbrace{\omega_{0,1}\otimes\dots\otimes\omega_{0,k_0}}_{\Omega_0}\otimes\cq{s_0}\\
    &\otimes\underbrace{\omega_{1,1}\otimes\dots\otimes\omega_{1,k_1}}_{\Omega_1}\otimes\cq{s_1}\\
    &\otimes\cdots\\
    &\otimes\underbrace{\omega_{n,1}\otimes\dots\otimes\omega_{n,k_n}}_{\Omega_n}\otimes\cq{s_n}%
\end{align*}
$s_i\in S$ and $\Omega_i$ range over the basis of $T_\ZZ(H^1(E,\CC))$ consisting of elements $\bigotimes_{j=1}^{k_i}\omega_{i,j}$ with $\omega_{i,j}\in\cq{\omega,\ol\omega}$. This generator lies in the component of $\CL_{E,S,s}^\vee$ of depth $n$ and weight $2n+\sum_{i=0}^nk_i$.

Suppose a tangent vector $v_s$ has been chosen at each $s\in S$. We assemble the $\CL_{E,S,v_s}^\vee$ as the base point $s$ ranges over $S$ into a Lie coalgebra
\[\widetilde\CL_{E,S}^\vee:=\bigoplus_{s\in S}\CL_{E,S,v_s}^\vee.\] All direct summands are isomorphic, but the maps $\Cor_\Hod$ on different components do not coincide. We will write $\Cor_s$ as a short notation for the map $\Cor_\Hod$ on the component corresponding to $s$,, extended so that $\Cor_s(s,\dots)=0$, i.e., the correlator of an element that contains the base point vanishes.

\subsubsection{Generating series}

We will package the correlators of depth $n$ into generating series in $2(n+1)$ commuting formal variables $t_0,\ol t_0,t_1,\ol t_1,\dots,t_n,\ol t_n$. We identify $t_i,\ol t_i$ with generators of $H_1(E,\ZZ)$ dual to $\omega,\ol\omega$. That is, the monomials in the $t_i,\ol t_i$ are identified with the generators of $\bigotimes_{i=0}^n\Sym(H_1(E,\ZZ))$.

For $x_0,\dots,x_n\in S$ and $s\in S$, define the generating series
\begin{align}
    &\Theta_s\pg{x_0:x_1:\dots:x_n}{t_0:t_1:\dots:t_n}
    =\,&\sum_{\Omega_0,\dots,\Omega_n}\Cor_s\pq{\Omega_0,x_0,\dots,\Omega_n,x_n)}(\Omega_0\du\otimes\dots\otimes\Omega_n\du),\label{eqn:def_theta}
\end{align}
where the sum is taken over the basis of $T_\ZZ(H^1(E,\CC))$ as above. The coefficient of $\prod_it_i^{m_i}\ol t_i^{m_i'}$ is the sum of all generators where $m_i$ copies of $\omega$ and $m_i'$ copies of $\ol\omega$ appear between $s_i$ and $s_{i+1}$. Letting $S_{m,m'}$ be the sum of generators of the degree-$(m,m')$ component of $T_\ZZ(H^1(E,\CC))$, i.e., the sum of all permutations of $\o^{\otimes m}\otimes\ob^{\otimes m'}$, this sum can be written
\begin{equation}
    \Cor_s\pq{S_{m_0,m_0'}\otimes(x_0)\otimes S_{m_1,m_1'}\otimes(x_1)\otimes\dots\otimes S_{m_n,m_n'}\otimes(x_n)}.
    \label{eqn:symm_corr}
\end{equation}
These coefficients are called the \emph{symmetric} Hodge correlators.

We also define, for $w_0,\dots,w_n\in E$ with $w_0+\dots+w_n=0_E$,
\begin{align*}
    \Theta\du_s\pg{w_0,w_1,\dots,w_n}{t_0:t_1:\dots:t_n}&=\Theta_s\pg{0:w_1:w_1+w_2:\dots:w_1+\dots+w_n}{t_0:t_1:\dots:t_n},
\end{align*}
and, for $u_0+\dots+u_n=0$,
\begin{align*}
    \Theta_s\pg{x_0:x_1:\dots:x_n}{u_0,u_1,\dots,u_n}&=\Theta\pg{x_0:x_1:\dots:x_n}{0:u_1:u_1+u_2:\dots,u_1+\dots+u_n}.
\end{align*}

The subspace generated by the elements of $\CLie_{E,S}^\vee$ having the form of the argument of (\ref{eqn:symm_corr}) is dual to a certain quotient of the Lie algebra $\Der^S(\gr^W\pi_1^{\rm nil}(E-S,v_0))$. This is the quotient by the image of the adjoint action of $H_1(E;\ZZ)$ mentioned in the introduction. In depth 0 and weight $>1$, the elements (\ref{eqn:symm_corr}) vanish, by the shuffle relations. In depth 0 and weight 1 -- i.e., elements $\Cor(\o_1,s_0)$ -- the elements are identified with elements $[s_0]-[s]$ in the Jacobian of $E$ (see \cite{goncharov-hodge-correlators}, \S10.5), and, in particular, vanish if $s$ and $s_0$ are torsion points. As we will see below, modulo the depth filtration, the symmetric ccorrelators form a subcoalgebra, as the terms $\delta_{\rm Cas}$ of the coproduct vanish.

Now let us establish some basic properties of the generating series.

\begin{lma}
\begin{enumerate}[(a)]
\item For $n>0$, the generating series $\Theta\pg{:}{:}$ are homogeneous in the $t_i$ and satisfy the dihedral symmetry relations:
\begin{align*}
&\Theta_s\pg{x_0:\dots:x_n}{t_0:\dots:t_n}\\
=\,&\Theta_s\pg{x_0+x:\dots:x_n+x}{t_0+t:\dots:t_n+t}&\text{(homogeneity)}\\
=\,&\Theta_s\pg{x_1:\dots:x_n:x_0}{t_1:\dots:t_n:t_0}&\text{(cyclic symmetry)}\\
=\,&(-1)^{n+1}\Theta_s\pg{x_n:\dots:x_1:x_0}{t_n:\dots:t_1:t_0}.&\text{(reflection)}
\end{align*}
\item For an automorphism $\phi\in\Aut(E)$,
\[
    \Theta_s\pg{x_0,\dots,x_n}{t_0:\dots:t_n}
    =
    \Theta_s\pg{\phi(x_0),\dots,\phi(x_n)}{\phi\cdot t_0:\dots:\phi\cdot t_n},
\]
where $\phi$ acts on the $t_i$ by the adjoint action on $H_1(E,\ZZ)$.
\item The elements $\Theta_s\pg{x_0:x_1:\dots:x_n}{u_0,u_1,\dots,u_n}$ satisfy the first shuffle relations:
\begin{equation}
    \sum_{\sigma\in\Sigma_{i,j}}\Theta_s\pg{x_{\sigma\inv(1)}:\dots:x_{\sigma\inv(i+j)}:x_0}{u_{\sigma\inv(1)},\dots,u_{\sigma\inv(i+j)},u_0}=0.
    \label{eqn:gf_first_shuffle}
\end{equation}
\end{enumerate}
\label{lma:theta_dihedral}
\end{lma}
\begin{proof}
The dihedral symmetry relations in (a) and the relation (b) are clear from the defnition of Hodge correlators.

The difficult part is homogeneity in $t_i$ and the first shuffle relation. For the former, it is enough to show
\[
\Theta_s\pg{x_0:\dots:x_n}{0:t_1:\dots:t_n}
=
\Theta_s\pg{x_0:\dots:x_n}{t_0:t_0+t_1:\dots:t_0+t_n}.
\]
Consider the coefficient of $\prod_it_i^{m_i'}\ol{t_i}^{m_i'}$ in the sum defining each side (\ref{eqn:def_theta}). For each $i$, fix an an ordering $\omega_{i,1}\dots\omega_{i,m_i+m_i'}$ of the word $\omega^{m_i}\ol\omega^{m_i'}$ and look at the terms in this coefficient in which the elements indexed by $t_i$ appear in the order specified by the word.

If $m_0=m_0'$, then both sides have exactly one such term \[\Cor_\Hod(x_0,1,x_1,\bigotimes_i\omega_{1,i}\dots,x_n,\bigotimes_i\omega_{n,i})),\] and they coincide. Otherwise, the coefficient on the left side is 0, while the terms on the right side are exactly the first shuffle relation on \[\Cor_\Hod(x_0,\underbrace{\bigotimes\omega_{0,i}},\underbrace{x_1,\bigotimes\omega_{1,i}\dots,x_n,\bigotimes\omega_{n,i})},\] which is 0. This proves homogeneity in the $t_i$.

Finally, (c) also follows from the first shuffle relation on the coefficients. To obtain the relation where $\cq{1,\dots,i}$ are shuffled with $\cq{i+1,\dots,i+j}$, we keep $x_0$ fixed and shuffle the $x_1,\dots,x_i$ and the forms indexed by $u_1,\dots,u_i$ with the other elements. (The proofs are identical for those for correlators on $\PP^1$; see \cite{malkin-shuffle}, Lemma 17.)
\end{proof}

\subsubsection{Coproduct}

The coproduct of the generating function $\Theta_s$ is in general difficult to write down. However, we can describe the terms of highest depth, which come from the $\delta_S$ component of the coproduct.
\begin{lma}
The coproduct of the generating functions $\Theta_s$ is given by
\begin{align*}
&\delta\Theta_s\pg{x_0:\dots:x_n}{t_0:\dots:t_n}=\\
= &\sum_{\rm cyc}\sum_{k=0}^n\Theta_s\pg{x_0:\dots:x_k}{t_0:\dots:t_k}\wedge\Theta_s\pg{x_k:x_{k+1}:\dots:x_n}{t_0:t_{k+1}:\dots:t_n} \\&+ \text{\rm lower depth terms}.
\end{align*}
The coproduct of the generating functions $\Theta_s\du$ is given by
\begin{align}
&\delta\Theta_s\du\pg{x_0,\dots,x_n}{t_0:\dots:t_n}=\nonumber\\
= &\sum_{\rm cyc}\sum_{k=0}^n\Theta_s\du\pg{-(x_1+\dots+x_k),x_1:\dots,x_k}{t_0:t_1:\dots:t_k}\wedge\Theta_s\du\pg{x_0,x_{k+1},\dots,x_n}{t_0:t_{k+1}:\dots:t_n} \nonumber\\&+ \text{\rm lower depth terms}.\label{eqn:thetastar_coproduct}
\end{align}
The lower-depth terms are Hodge correlators of elements that do not depend on $s$.
\end{lma}
\begin{proof}
The formula for the coproduct of $\Theta_s$ arises from the definition of the $\delta_S$ term of the coproduct. The formula for the coproduct of $\Theta_s\du$ would follow immediately from that for $\Theta_s$ if the $\Theta_s$ were invariant under an additive shift of the arguments $x_i$. This is Theorem~\ref{thm:depth_indep} below, which is independent of  (\ref{eqn:thetastar_coproduct}).
\end{proof}
These formulas for the coproduct formally coincide with those for the dihedral Lie coalgebra, defined by A.Goncharov in \cite{goncharov-dihedral} in order to study multiple polylogarithms, as well as in the quasidihedral Lie coalgebra modulo the depth filtration, defined by the author in \cite{malkin-shuffle} to study Hodge correlators on $\PP^1$.

\subsection{Symmetric correlators modulo depth}
\label{sec:mc_modulo_depth}

In this section, $H^1(X)$ always refers to $H^1(X;\CC)$.

\subsubsection{Change of base point formula}

Fix $p\in S$. Let us define a map $\rho_p:T(H^1(X))\to T(H^1(X)\oplus \QQ[S])$ as follows.

For a word $\omega_1\otimes\dots\otimes\omega_n\in T(H^1(X))$,
\[\rho_p(\omega_1\otimes\dots\otimes\omega_n)=\sum_k(-1)^k\sum_{\stackrel{i_1<i_2,\dots<i_k<n}{i_{j+1}>i_j+1}}\omega_1\otimes\dots\otimes(\gen{\omega_{i_j},\omega_{i_j+1}}(p))\otimes\dots\otimes\omega_n,\] where $\gen{,}$ is the skew-symmetric pairing: $\gen{\o,\ob}=-\gen{\ob,\o}=1$. That is, it is the sum over all possible replacements of pairs $(\o\otimes\ob)$ and $(\ob\otimes\o)$ by the puncture $p$, taken with appropriate sign. For example, we have:
\begin{align*}
\rho_p(1)&=1,\\
\rho_p(\o)&=\o,\\
\rho_p(\o\otimes\o)&=\o\otimes\o,\\
\rho_p(\o\otimes\ob)&=(\o\otimes\ob)-(p),\\
\rho_p(\o\otimes\ob\otimes\o)&=(\o\otimes\ob\otimes\o)-(p\otimes\o)+(\o\otimes p),\\
\rho_p(\o\otimes\ob\otimes\o\otimes\ob)&=(\o\otimes\ob\otimes\o\otimes\ob)-(p\otimes\o\otimes\ob)-(\o\otimes\ob\otimes p)+(\o\otimes p\otimes\ob)+(p\otimes p).
\end{align*}
For $a\in S$, define $\rho_p(a)=(a)-(p)$, extended by linearity to $\QQ[S]$. Then, extend $\rho_p$ to $CT(H^1(X)\oplus\QQ[S])$: if  $x_0,\dots,x_k\in\QQ[S]$, and $\Omega_0,\dots,\Omega_k\in T(H^1(X))$, then 
\[\rho_p(\Omega_0\otimes x_0\otimes\dots\otimes \Omega_k\otimes x_k)=\rho_p(\Omega_0)\otimes\rho_p(x_0)\otimes\dots\otimes\rho_p(\Omega_k)\otimes\rho_p(x_k).\]
\begin{lma}[Change of base point formula]
Suppose that $p\neq q$. Then the following relation holds for Hodge correlators in weight $>2$:
\begin{equation}
\Cor_p(x)=\Cor_q\pq{\rho_p(x)}.
\label{eqn:change_base_point}
\end{equation}
\label{lma:change_base_point}
\end{lma}
On the right side stands a sum of correlators obtained from the one on the left by taking all possible replacements of punctures and pairs of adjacent cohomology classes with $(p)$, taken with the appropriate sign.

Before proceeding to the proof, let us illustrate the formula on some examples. In weight $4$,
\begin{align*}
\Cor_p(a,b,c)&=\Cor_q(a,b,c)-\Cor_q(p,b,c)-\Cor_q(a,p,c)-\Cor_q(a,b,p),\\
\Cor_p(a,b,\o,\ob)&=\Cor_q(a,b,\o,\ob)\\&\quad-\Cor_q(p,b,\o,\ob)-\Cor_q(a,p,\o,\ob)-\Cor_q(a,b,p)+\Cor_q(p,p,\o,\ob),\\
\Cor_p(a,\o,\ob,\o,\ob)&=\Cor_q(a,\o,\ob,\o,\ob)\\&\quad-\Cor_q(p,\o,\ob,\o,\ob)-\Cor_q(a,p,\o,\ob)+\Cor_q(a,\o,p,\ob)-\Cor_q(a,\o,\ob,p)\\&\quad+\Cor_q(p,p,\o,\ob)-\Cor_q(p,\o,p,\ob)+\Cor_q(p,\o,\ob,p).
\end{align*}
If the left side of the expression only contains punctures, we recover a formula identical to the one found by \cite{goncharov-rudenko}, Theorem 2.6, for Hodge correlators on the punctured $\PP^1$. More generally, for symmetric correlators, we have:
\begin{cly}
Suppose that $p\neq q$. Then we have the relation in weight $>2$:
\begin{align*}
&\Cor_p(S_{m_0,m_0'}\otimes x_0\otimes S_{m_1,m_1'}\otimes x_1\otimes\dots\otimes S_{m_n,m_n'}\otimes x_n)
\\&=
\Cor_q(S_{m_0,m_0'}\otimes((x_0)-(p))\otimes\dots\otimes S_{m_n,m_n'}((x_n)-(p)))
\\&=
\sum_k(-1)^k\sum_{i_1<\dots<i_k}\Cor_q(S_{m_0,m_0'}\otimes x_0\otimes \dots \otimes p\otimes\dots\otimes p \otimes\dots\otimes S_{m_n,m_n'}\otimes x_n),
\end{align*}
where on the right the punctures $x_{i_1},\dots,x_{i_k}$ are replaced with $q$.
\label{cly:sym_change_base_point}
\end{cly}
\begin{proof}
For all $m,m'\geq0$, $\rho_p(S_{m,m'})=S_{m,m'}$.
\end{proof}

\begin{proof}[Proof of Lemma~\ref{lma:change_base_point}]
We first prove the change of base point formula in the real Hodge realization, i.e, that it holds on the level of the Hodge correlator functions $\Cor_\HHH$.

The Green's functions associated to the points $p$ and $q$ are related by \[G_p(x,y)=G_q(x,y)-G_q(x,p)-G_q(y,p)+C,\] where $C$ is a constant that depends on the choices of tangent vectors at $p$ and $q$. Now consider any tree contributing to the Hodge correlator of $\Omega_0\otimes(x_0)\otimes\dots\otimes\Omega_k\otimes x_k$. Write the Green's function $G_p(x,y)$ assigned to each edge in terms of the $G_q$, and examine the contribution of the three terms in $G_p(x,y)-G_q(x,y)$: $G_q(x,p)$, $G_q(y,p)$, and $C$ for a given edge. There are three cases:
\begin{enumerate}[(1)]
\item An external edge $E$ decorated by a puncture $a$, assigned the function $G_p(x,a)$. Assigning the form $-G_q(x,p)$ to $E$ gives the correlator where $a$ has been replaced by $-(p)$. The terms $C$ and $G_q(a,p)$ are constants. Because the Hodge correlator has weight $>2$, there is at least one internal edge in the tree, so the correlator where a constant has been placed on $E$ is the integral of an exact form $d^\CC(\dots)$.
\item An internal edge $E$ that splits the tree into two parts, one of which is decorated by two 1-forms. Suppose that in $G_p(x,y)$, the vertex assigned the variable $x$ is adjacent to external vertices labeled $\omega_1$ and $\omega_2$. Then the terms $G_q(x,p)$ and $C$ are independent of $y$, and the integral splits into a product; the integrand for the subtree growing from $y$ is an exact form, so we get 0. For the term $-G_q(y,p)$, the integral also splits into a product of $\int_E\omega_1\wedge\omega_2$ and the correlator with $x$ replaced by an external vertex $-(p)$.
\item An internal edge $E$ that splits the tree into two parts, each of which is decorated by at least one puncture. Then, as in the previous case, each term in the expression for $G_p(x,y)$ is independent of either $x$ or $y$. The integral splits into a product of two factors, one of which is 0.
\end{enumerate}
We conclude that the change of base point is computed by adding all possible replacements of external punctures $a$ by $-(p)$ and pairs $\omega_1\otimes\omega_2$ by $-\gen{\omega_1,\omega_2}(p)$. This implies the lemma.

(Note that the assumption of weight $>2$ was crucial to all arguments involving integration of the exact form.)

One easily checks by induction that the coproducts of the two sides of (\ref{eqn:change_base_point}) are equal. This implies the result on the level of the Hodge correlators $\Cor_\Hod$.
\end{proof}

\subsubsection{Independence on base point}

In this part, we prove the following important result.
\begin{thm}
The symmetric Hodge correlators in weight $>2$ are independent of the base point modulo the depth filtration.
    
Precisely, let $x\in CT(H^1(X)\oplus\QQ[S])$. Then there exists $\widetilde x$, equal to $x$ modulo lower-depth terms, such that $\Cor_p(\widetilde x)$ is independent of $p$.
\label{thm:depth_indep}
\end{thm}

In terms of generating functions, this theorem implies:
\begin{cly}
The generating functions $\Theta\du$ satisfy the dihedral symmetry relations of Lemma~\ref{lma:theta_dihedral}:
\begin{align*}
    &\Theta\du_s\pg{w_0,\dots,w_n}{t_0:\dots:t_n}\\
    =\,&\Theta\du_s\pg{w_1,\dots,w_n,w_0}{t_1:\dots:t_n:t_0}\\
    =\,&(-1)^{n+1}\Theta\du_s\pg{w_n,\dots,w_1,w_0}{t_n:\dots:t_1:t_0}
\end{align*}
modulo lower-depth terms that are independent of $s$.
\label{cly:thetastar_dihedral}
\end{cly}
\begin{proof}
By the cyclic symmetry and dihedral relations on correlators, these expressions are equal up to an additive shift in the correlators' arguments, equivalently, a change in base point.
\end{proof}

Notice that all terms on the right side of (\ref{eqn:change_base_point}) have \emph{higher}  or equal depth to the left side. It will be necessary to find correction terms of lower depth to obtain a formula of the form \[\Cor_p(h_0\otimes x_0\otimes\dots\otimes h_k\otimes x_k)+\Cor_p(\text{lower depth})=\Cor_q(h_0\otimes x_0\otimes\dots\otimes h_k\otimes x_k)+\Cor_q(\text{l.d.})\] when each $h_i$ a symmetric expression $S_{m,m'}$.

The proof of the theorem relies on a key construction. We will find elements:
\[S_{m_0,m_0'}*S_{m_1,m_1'}*\dots*S_{m_n,m_n'}\in T(H^1(X))\]
such that
\begin{align}
&\rho_p(S_{m_0,m_0'}*S_{m_1,m_1'}*\dots*S_{m_n,m_n'})=\label{eqn:star_property}
\\&=
\sum_k\sum_{i_1<\dots<i_k} \pq{S_{m_0,m_0'}*\dots*S_{m_{i_1-1},m_{i_1-1}'}}\otimes(p)\otimes\pq{S_{m_{i_1},m_{i_1}'}*\dots*S_{m_{i_2-1},m_{i_2-1}'}}\otimes(p)\otimes\dots.\nonumber
\end{align}
Before showing how to construct these elements, let us prove the theorem, assuming these elements exist.

\begin{proof}[Proof of Theorem~\ref{thm:depth_indep}]
Consider an element 
\[x=S_{m_0,m_0'}\otimes x_0\otimes S_{m_1,m_1'}\otimes x_1\otimes\dots\otimes S_{m_n,m_n'}\otimes x_n.\]
Let $I$ be a proper subset of $\cq{0,\dots,n}$. Write $I$ as the union of its cyclically contiguous subsets, each of the form $\cq{i,i+1,\dots,i+k}$ (indices modulo $n+1$). Let $x_{/I}$ be the element formed by replacing each \[S_{m_i,m_i'}\otimes x_i\otimes\dots\otimes x_{i+k}\otimes S_{m_{i+k},m_{i+k}'}\] by $S_{m_i,m_i'}*\dots*S_{m_{i+k},m_{i+k}'}$.

Now consider the corrected element:
\[\widetilde x=\sum_I(-1)^{\aq I}x_{/I}.\]
It is equal to $x$ modulo the depth filtration. Also, let \[y_q=S_{m_0,m_0'}\otimes q\otimes S_{m_1,m_1'}\otimes q\otimes\dots\otimes S_{m_n,m_n'}\otimes q,\] and define $\widetilde y_q$ in the same way. By a standard inclusion-exclusion argument, the property (\ref{eqn:star_property}) implies that
\[\rho_p(\widetilde x+\widetilde y_q)=\widetilde x + \text{(terms containing $q$)}.\]
Because the correlator with base point $q$ is zero for the terms containing $q$, this gives
\[\Cor_p(\widetilde x+\widetilde y_q) = \Cor_q(\widetilde x).\] On the other hand, the Hodge correlator $\Cor_p(\widetilde y_q)$ depends only on $p-q$, and thus $p\mapsto\Cor_p(\widetilde x)-\Cor_0(\widetilde x)$ provides a group homomorphism $E\to\RR$, and must be 0. Therefore, $\Cor_p(\widetilde x)$ is independent of $p$.
\end{proof}
    
\begin{lma}
There exist elements, independent of choice of symplectic basis of $H^1(X)$, satisfying (\ref{eqn:star_property}).
\label{lma:depth_correction_term}
\end{lma}
\begin{proof}
We produce such elements explicitly:
\[S_{m_0,m_0'}*\dots*S_{m_k,m_k'}=\f{1}{2^k}\sum_{n_0,n_0',\dots,n_k,n_k'}\pm S_{n_0,n_0'}\otimes S_{n_1,n_1'}\otimes\dots\otimes S_{n_k,n_k'},\]
where the sum is taken over the $n_i,n_i'\geq0$ such that:
\begin{align*}
n_i+n_i'&=\begin{cases}m_i+m_i'+1&i=0,k\\m_i+m_i'+2&0<i<k\end{cases},\\
(n_0-n_0')+\dots+(n_k-n_k')&=(m_0-m_0')+\dots+(m_k-m_k').
\end{align*}
A term is taken with the sign $-$ if there is an odd number of $i$ ($i=0,\dots,k-1$) such that
\[(n_0-n_0')+\dots+(n_i-n_i')<(m_0-m_0')+\dots+(m_i-m_i'),\]
otherwise with the sign $+$.

Examples:
\begin{align*}
S_{0,0}*S_{0,0}&=\f12\pq{\o\ob-\ob\o},\\
S_{0,0}*S_{1,0}&=\f12\pq{\o\o\ob+\o\ob\o-\ob\o\o},\\
S_{0,0}*S_{0,0}*S_{0,0}&=\f14\pq{\o\o\ob\ob+\o\ob\o\ob-\o\ob\ob\o-\ob\o\o\ob+\ob\o\ob\o+\ob\ob\o\o}.
\end{align*}

We explain the construction by picture. The basis elements of $T(H^1(X))$ of a given weight are in bijection with lattice paths: a word $\o_1\otimes\dots\otimes\omega_n$ corresponds to the path whose $i$-th step is $(1,0)$ if $\o_i=\o$ and $(0,1)$ if $\o_i=\ob$. The elements of $T(H^1(X)\oplus\QQ[p])$ are lattice paths that also allow the diagonal step $(1,1)$, corresponding to $(p)$. (The points of the lattice path are simply the Hodge bidegrees of the initial subwords.) The map $\rho_p$ replaces a path by the sum of all paths obtained by replacing steps (up, right) or (right, up) with diagonal steps, in the latter case changing the sign. 

\begin{figure}[t]
\begin{tabular}{cc}
\includegraphics[height=0.25\textwidth]{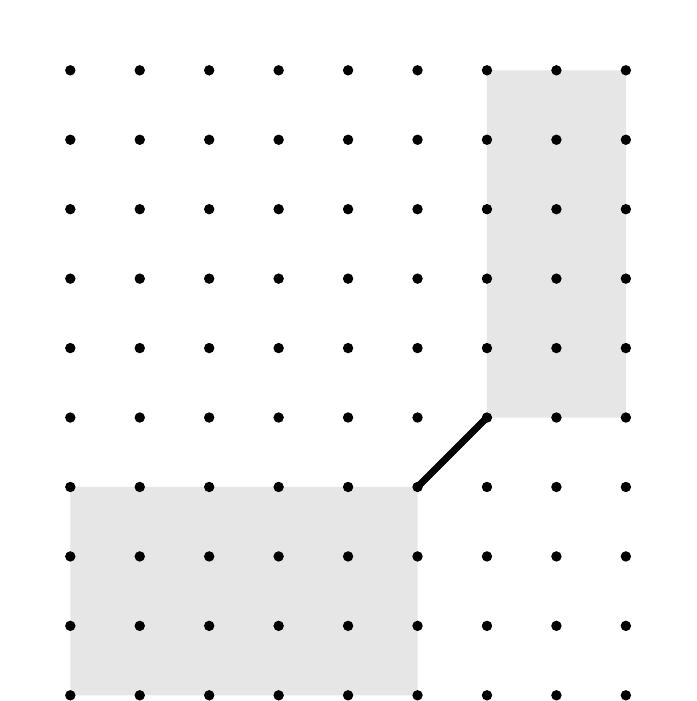}&
\includegraphics[height=0.25\textwidth]{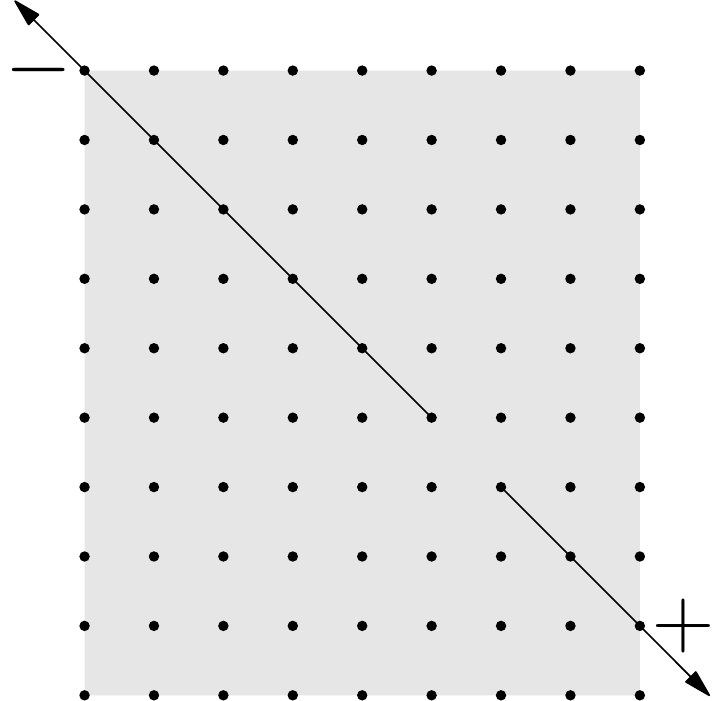}\\
$S_{m_0,m_0'}\otimes(p)\otimes S_{m_1,m_1'}$&$S_{m_0,m_0'} * S_{m_1,m_1'}$
\end{tabular}
\simplecap{fig:lp1}{Construction of the element $S_{m_0,m_0'}*S_{m_1,m_1'}$: the paths crossing the rays marked $-$ and $+$ are taken with the corresponding sign.}
\end{figure}

To construct the element, we first consider the concatentation of paths in $S_{m_0,m_0'}$, \dots, $S_{m_k,m_k'}$, with a step $(1,1)$ inserted between each pair. Draw a diagonal line $\ell_i$ bisecting the step that was inserted between $S_{m_{i-1},m_{i-1}'}$ and $S_{m_i,m_i'}$. Any path of the Hodge bidegree $\pq{\sum m_i+k,\sum m_i'+k}$ appears in a unique term $S_{n_0,n_0'}\otimes\dots\otimes S_{n_k,n_k'}$, and each $S_{n_i,n_i'}$ is the sum of paths between the lines $\ell_i$ and $\ell_{i+1}$. The sign of a path is determined by the rays on which it crosses the diagonal lines: $+$ if below the step, $-$ if above. See Figure~\ref{fig:lp1}.

\begin{figure}[t]
\includegraphics[height=0.25\textwidth]{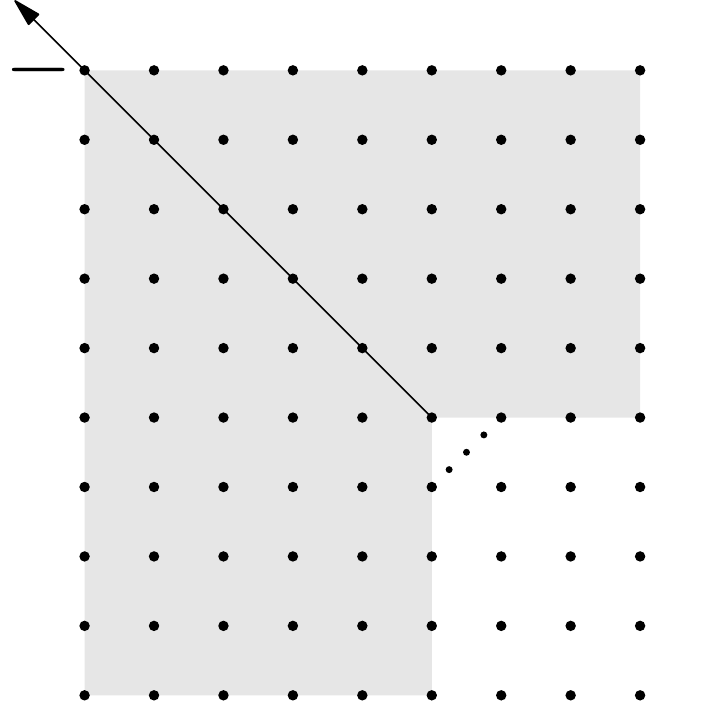}
\includegraphics[height=0.25\textwidth]{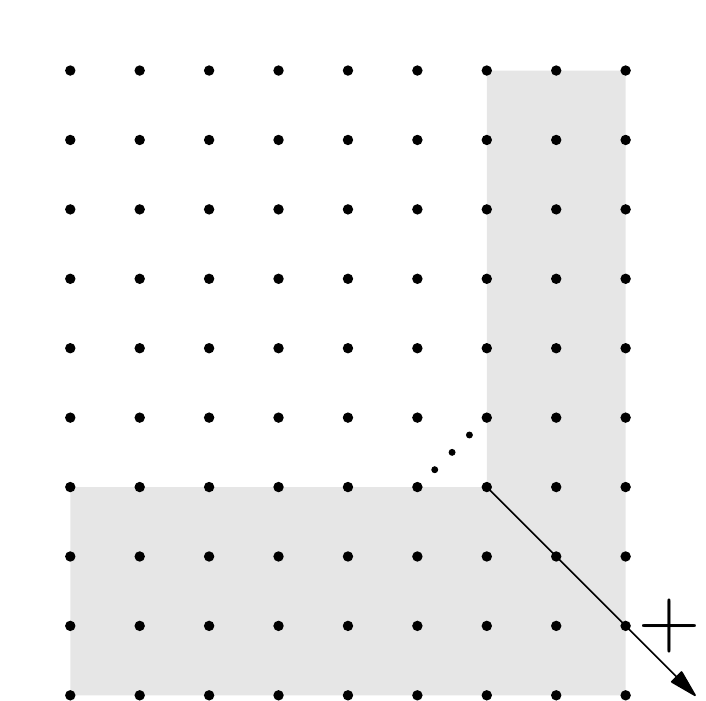}
\simplecap{fig:lp2}{The point of nonconcavity contributing a term to the right side of (\ref{eqn:star_property}).}
\end{figure}

Now fix a choice of a ray of each such diagonal, and consider the terms coming from lattice paths crossing these rays. We claim that any such term satisfies (\ref{eqn:star_property}) modified by a factor of $\f{1}{2^k}$. Indeed, these are the lattice paths lying in a certain rectilinear region (right part of the figure). Most terms in $\rho_p$ are canceled; the only terms remaining are those with segments $(1,1)$ at the points of nonconcavity of this region. This is precisely the expression on the right of (\ref{eqn:star_property}). See Figure~\ref{fig:lp2}.
\end{proof}

The simplest example of the corrected correlator, for $(a)\otimes(b)\otimes(c)$:
\begin{align*}
(a)\otimes(b)\otimes(c)&-\f12\pq{\o\otimes\ob-\ob\otimes\o}\otimes(b)\otimes(c)\\&-\f12(a)\otimes\pq{\o\otimes\ob-\ob\otimes\o}\otimes(c)\\&-\f12(a)\otimes(b)\otimes\pq{\o\otimes\ob-\ob\otimes\o}.
\end{align*}
(The terms where two points were replaced are 0, because of the reflection relations.)

\subsection{Second shuffle relations}

\subsubsection{The depth 2 case: dihedral symmetry}

\begin{thm}
The corrected symmetric Hodge correlators in depth 2 satisfy the second shuffle (dihedral symmetry) relations modulo terms of lower depth that are independent of the base point.

Precisely, the corrected element for
\begin{align}
    &S_{m_0,m_0'}\otimes(0)\otimes S_{m_1,m_1'}\otimes(x_1)\otimes S_{m_2,m_2'}\otimes(x_1+x_2) \nonumber\\+\, &S_{m_0,m_0'}\otimes(0)\otimes S_{m_2,m_2'}\otimes(x_2)\otimes S_{m_1,m_1'}\otimes(x_1+x_2)\label{eqn:elliptic_second_shuffle}
\end{align}
lies in the kernel of the map $\Cor_s$ for every $s$.
\label{thm:dihedral_depth2}
\end{thm}
\begin{proof}
The corrected element for (\ref{eqn:elliptic_second_shuffle} changes sign under the map $x\mapsto(x_1+x_2-x)$ and reflection. On the other hand, it is invariant under this operation up to an additive shift (i.e., change in base point).
\end{proof}

\subsubsection{Relations in higher depth}

The \emph{second shuffle relations} are relations of the form
\[\sum_{\sigma\in\Sigma_{i,j}}\Theta\du\pg{x_0,x_{\sigma\inv(i)},x_{\sigma\inv(2)},\dots,x_{\sigma\inv(i+j)})}{t_0,t_{\sigma\inv(1)},t_{\sigma\inv(2)},\dots,t_{\sigma\inv(i+j)}}+\dots,\]
perhaps with additional terms of lower depth. The Hodge correlators on $\PP^1$ are known to obey such relations, in addition to the first shuffle relations, the structural relations in $\CL^\vee_{X,S,v_0}$; the lower-depth terms were described precisely by \cite{malkin-shuffle}.

The relation of Theorem~\ref{thm:dihedral_depth2} is a special case of a second shuffle relation. In depth $>2$, the second shuffle relations are not equivalent to dihedral symmetry. However, one hopes for a generalization. 
\begin{conj}
The second shuffle relations for symmetric elliptic Hodge correlators hold modulo the depth filtration. The lower-depth terms are independent of the base point $s$.
\label{cnj:second_shuffle}
\end{conj}
The lower-depth correction terms in depth $>2$ are not known. In particular, the corrected correlators do not satisfy the second shuffle relations in higher depth. However, calculations in low weight support this conjecture. We may expect the elliptic relations to be deformations of the relations for $\PP^1$ (see \S\ref{sec:degen}).

\section{Bianchi hyperbolic threefolds and modular complexes}

\label{sec:bianchi_and_modular}

\subsection{Bianchi tesselations and orbifolds}

\subsubsection{Definition}

Let $K=\QQ[\sqrt{-d}]$ be an imaginary quadratic field with lattice of integers $\OOO$. The \emph{Bianchi tesselation} (\cite{bianchi}) is an ideal polyhedral tesselation of the upper half-space $\HH^3$ associated with $\OOO$, whose cell complex has a natural structure of a complex of $\GL_2(\OOO)$-modules. We define it now.

Let $\ol\FFF$ be the space of positive semidefinite Hermitian forms on $(\OOO^2\otimes_\OOO\CC)^*$. The subset $\FFF$ of positive definite forms is a dense open subset of $\ol\FFF$. We identify $\HH^3$ and its compactification $\ol\HH^3=\HH^3\cup\PP^1(\CC)$ with the real projectivizations of $\FFF$ and $\ol\FFF$, respectively. The action of $\GL_2(\CC)$ on $\CC^2$ provides an action on $\ol\FFF$ that descends to an action on $\ol\HH^3$.

Every $v\in\OOO^2$ provides a positive semidefinite form $\aq{\gen{-,v}}^2\in\del\ol\FFF$. The convex hull of the set \[\cq{\aq{\gen{-,v}}^2:\text{$v$ a primitive vector in $\OOO^2$}}\] is a polyhedron in $\ol\FFF$ with vertices on the boundary. The polyhedron projects to an ideal tesselation of $\HH^3$ with vertices on $\PP^1(\OOO)\subset\PP^1(\CC)$. Let $B\chain$ be the polyhedral cell complex over $\ZZ$ of this ideal tesselation. We will shift this complex in degree so that the space of $i$-dimensional cells it in degree $3-i$ ($i=0,1,2,3$). We get a cohomological complex \[B^0\tto\del B^1\tto\del B^2\tto\del B^3.\]
The group $\GL_2(\OOO)$ acts on the Bianchi tesselation, giving $B\chain$ the structure of a complex of left $\GL_2(\OOO)$ modules.

The quotient $\GL_2(\OOO)\sm\HH^3$ is a finite-volume hyperbolic threefold with cusps in bijection with the ideal class group of $\OOO$. If $\Gamma$ is a finite-index subgroup of $\GL_2(\OOO)$, the quotient $\Gamma\sm\HH^3$ is also a finite-volume hyperbolic threefold with a finite map to $\GL_2(\OOO)\sm\HH^3$. 

A right $\GL_2(\OOO)$-module $T$ provides a local system on $\Gamma\sm\HH^3$, which we also denote by $T$. Then the chain complex of $\GL_2(\OOO)\sm\HH^3$ with coefficients in $T$ is \begin{equation}T\otimes_\Gamma B\chain\cong\pq{\ZZ[\Gamma\sm\GL_2(\OOO)]\otimes T}\otimes_{\GL_2(\OOO)}B\chain.\label{eqn:loc_sys_cc}\end{equation}

\subsubsection{The Gaussian and Eisenstein cases}

Following \cite{goncharov-euler}, for $d=1$ ($\OOO=\ZZ[i]$) and $d=3$ ($\OOO=\ZZ[\rho]$) we have the following description of the Bianchi complexes in degrees 1 and 2. 

The action of $\GL_2(\OOO)$ is transitive on the $i$-dimensional cells for each $i$. Choose $\GL_2(\OOO)$-generators $G_i\in B^i$: we may take
\begin{align*}
G_1&=\text{(the ideal triangle $(1,0,\infty)$)}\\
G_2&=\text{(the geodesic $(0,\infty)$)}
\end{align*}
where $(v_1,\dots,v_n)$, $v_i\in\PP^1(\OOO)=\PP(V^2(\OOO))$, denotes the oriented cell with ideal vertices at $v_1,\dots,v_n$ under the identification of $\PP^1(\CC)$ with the boundary of $\ol\HH^3$.
Let $D_i$ be the subgroup of $\GL_2(\OOO)$ stabilizing $G_i$. 

The group $D_1$ stabilizing the triangle $(0,1,\infty)$ is isomorphic to \[S_3\times \OOO\uns.\] The first component $S_3$ acts on $(v,w)\in \OOO\oplus\OOO$ by permutations of the triple $(v,w,-v-w)$, i.e., the generators of $S_3$ are represented by
\[(123)\mapsto\begin{pmatrix}0&-1\\1&-1\end{pmatrix},\quad(12)\mapsto\begin{pmatrix}0&1\\1&0\end{pmatrix}.\] The second component acts by scalars. There is a sign homomorphism $\chi_1:D_1\to\ZZ$ keeping track of the action of $D_1$ on the orientation, with $\chi_1((123))=1$ and $\chi_1((12))=-1$. So the space of 2-cells is \[B^1=\ZZ[GL_2(\OOO)]\otimes_{D_1}\chi_1.\]

The group $D_2$ stabilizing the geodesic $(0,\infty)$ is isomorphic to \[S_2\ltimes(\OOO\uns\times\OOO\uns),\] with $S_2$ acting on $\OOO\uns\times\OOO\uns$ by permutation of the factors. The nontrivial element of $S_2$ acts by $(v,w)\mapsto(w,v)$ and $\OOO\uns\times\OOO\uns$ acts diagonally. There is a sign homomorphism $\chi_2:D_2\to\ZZ$, and the space of 1-cells is \[B^2=\ZZ[GL_2(\ZZ[i])]\otimes_{D_2}\chi_2.\]

Let $\Ppp$ be a prime ideal in $\OOO$. The group $\GL_2(\OOO)$ acts on the quotient $(\ZZ[i]/\Ppp)^2$. Let $\Gamma_1(\Ppp)$ be the stabilizer in $\GL_2(\OOO)$ of the vector $(0,1)\in(\ZZ[i]/\Ppp)^2$. The action on the vector $(0,1)$ provides an isomorphism of $\GL_2(\OOO)$-modules \[\ZZ[\Gamma_1(\Ppp)\sm\GL_2(\OOO)]\cong\ZZ[\FF_\Ppp^2-0],\quad\FF_\Ppp=\OOO/\Ppp.\]

The chain complex (\ref{eqn:loc_sys_cc}) of $\Gamma_1(\Ppp)\sm\HH^3$ with coefficients in a local system $T$ is then identified in degrees 1 and 2 with 
\begin{align*}
T\otimes_{\Gamma_1(\Ppp)}B\chain
&\cong\pq{\ZZ[\Gamma_1(\Ppp)\sm\GL_2(\OOO)]\otimes T}\otimes_{\GL_2(\OOO)}B\chain
\\
&\cong\pq{\ZZ[\FF_\Ppp^2-0]\otimes T}\otimes_{\GL_2(\OOO)}\pq{\ZZ[\GL_2(\OOO)]\otimes_{D_\bullet}\chi_\bullet}.%
\end{align*}
This space is generated in degree $i$ by elements 
\[\pq{(\alpha,\beta)\otimes t}\otimes(G_i),\quad(\alpha,\beta)\in\FF_\Ppp^2-0,\quad t\in T.\]

\subsection{Modular complexes}

\subsubsection{Definition}

Let $\OOO=\ZZ$ or the lattice of integers in an imaginary quadratic field. We are going to define the modular complexes $M\chain_k$, complexes of left $\GL_k(\OOO)$-modules that generalize the complexes defined by \cite{goncharov-polylogs-modular} for $\GL_k(\ZZ)$.

Fix a $k$-dimensional $\OOO$-vector space $V$. An \emph{extended basis} of $V$ is a sequence of vectors $\gen{v_0,v_1,\dots,v_k}$, $v_i\in V$, such that $v_0+\dots+v_k=0$ and $v_1,\dots,v_k$ form a basis of $V$. (Consequently, any other set of $k$ vectors in this sequence form a basis.) We also use the notation
\begin{align*}
\bq{v_1,\dots,v_k}&=\gen{-v_1-\dots-v_k,v_1,v_2,\dots,v_k},\\
\bq{v_1:\dots:v_k}&=\bq{v_2-v_1,v_3-v_2,\dots,v_k-v_{k-1},-v_k}.
\end{align*}
The set $B_V$ of extended bases of $V$ is a principal homogeneous space for $\GL(V)$.

The complex of left $\GL_k(\OOO)$-modules $M\chain_k$ lies in the degrees $1,\dots,n$. The module $M^1_k$ is the quotient of $\ZZ[B_V]$ by the double shuffle relations
\begin{align}
    \sum_{\sigma\in\Sigma_{i,j}}\bq{v_{\sigma\inv(1)}:\dots:v_{\sigma\inv(i+j)}}&=0,&\text{\it(first shuffle)} \label{eqn:modular_first_shuffle}\\
    \sum_{\sigma\in\Sigma_{i,j}}\bq{v_{\sigma\inv(1)},\dots,v_{\sigma\inv(i+j)}}&=0.&\text{\it(second shuffle)}\label{eqn:modular_second_shuffle}
\end{align}
\begin{lma}[\cite{goncharov-polylogs-modular}, Theorem 4.1]
    The double shuffle relations imply the dihedral symmetry relations:
    \begin{align}
        \gen{v_0,v_1,\dots,v_k}=\gen{v_1,\dots,v_k,v_0}=(-1)^{k+1}\gen{v_k,\dots,v_1,v_0}=\gen{-v_0,-v_1,\dots,-v_k}.\label{eqn:modular_dihedral}
    \end{align}
\end{lma}

The module $M^n_k$ is generated by elements
\[[v_1,\dots,v_{k_1}]\wedge\dots\wedge[v_{k_{n-1}+1},\dots,v_{k_n}],\]
where each block $\bq{v_{k_{i-1}+1},\dots,v_{k_i}}$ is an extended basis of a sublattice $V_i$ in $V$, and $V=V_1\oplus\dots\oplus V_n$ (from which is follows that $k_1+\dots+k_n=k$). The double shuffle relations are imposed on each of the blocks, and the blocks anticommute.

The coproduct $\delta:M^1_k\to M^2_k$ is defined by 
\[\delta\gen{v_0,v_1,\dots,v_k}=\sum_{\rm cyc}\sum_{i=1}^k[v_0,\dots,v_{i-1}]\wedge[v_{i+1},\dots,v_k]\]
with the outer cyclic sum is over $\cq{0,1,\dots,k}$. The coproduct is extended by the Leibniz rule to the higher degrees, i.e., 
\[\delta(x_1\wedge\dots\wedge x_n)=\sum_{i=1}^n(-1)^{i+1}x_1\wedge\dots\wedge\delta(x_i)\wedge\dots\wedge x_n.\]

We will also consider the \emph{relaxed modular complex} $\widetilde M_k^n$, in which impose only the first shuffle relations (\ref{eqn:modular_first_shuffle}) and the dihedral symmetry relations (\ref{eqn:modular_dihedral}). By the lemma, the modular complex is the quotient of the relaxed modular complex by the second shuffle relations (\ref{eqn:modular_second_shuffle}).

\subsubsection{Relating the Gaussian and Eisenstein Bianchi and modular complexes for $k=2$}

In this section, suppose $\OOO=\ZZ[i]$ or $\ZZ[\rho]$. We will construct an isomorphism between the modular complex $M_2\ch$ and the Bianchi complex $B\ch$ in degrees 1 and 2.

Recall that $B\ch$ is generated by the ideal triangle $(1,0,\infty)$ in degree 1 and the geodesic $(0,\infty)$ in degree 2, with the boundary map given by 
\[(1,0,\infty)\mapsto(1,0)+(0,\infty)+(\infty,1).\]
The modular complex $M\ch_2$ is generated in degree 1 by the extended basis $[e_1,e_2]$, with the coproduct
\[\bq{e_1,e_2}\mapsto\bq{-e_1-e_2}\wedge\bq{e_2}+\bq{e_1}\wedge\bq{-e_1-e_2}+\bq{e_2}\wedge\bq{e_1}.\]
Making as before the identification of $\PP^1(\OOO)$ with $\PP^1(V)$, define the map $\psi:M\ch_2\to B\ch$ by
\[
    \psi\pq{\gen{v_1,v_2,v_3}}=\text{the triangle $(v_1,v_2,v_3)$},\quad \psi\pq{[v_1]\wedge[v_2]}=\text{the geodesic $(v_1,v_2)$}.
\]
\begin{lma}
The map $\psi$ is an isomorphism of complexes of $\GL_2(\OOO)$-modules.
\label{lma:modular_bianchi}
\end{lma}
\begin{proof}
By construction, $\psi$ is a surjective map of abelian groups. We must verify (1) $\psi$ commutes with the action of $\GL_2(\OOO)$, (2) $\psi$ commutes with the coproduct, (3) $\psi$ respects the double shuffle relations, and the images of the double shuffle and anticommutation relations are all relations in $B\ch$.

(1) holds by construction. For (2), notice that
\begin{align*}
    \delta\bq{e_1,e_2}
    &=\bq{-e_1-e_2}\wedge\bq{e_2}+\bq{e_1}\wedge\bq{-e_1-e_2}+\bq{e_2}\wedge\bq{e_1}.\\
    &=\bq{e_2}\wedge\bq{e_1} + \begin{pmatrix}0&-1\\1&-1\end{pmatrix}\bq{e_2}\wedge\bq{e_1} + \begin{pmatrix}0&-1\\1&-1\end{pmatrix}^2\bq{e_2}\wedge\bq{e_1}
\end{align*}
and that $\begin{pmatrix}0&-1\\1&-1\end{pmatrix}$ acts by cyclic permutation on $(0,1,\infty)$. 
    
For (3), double shuffle relation in $M^1_2$ is just equivalent to dihedral symmetry, which is precisely the relation imposed by $\otimes_{D_1}\chi_1$. The only relations in $M_2^2$ are the anticommutation relation and the relation $[v_1]=[-v_1]$, whose images are the only relations among the 1-cells in $B^2$.
\end{proof}

As a consequence, the chain complex of $\Gamma_1(\Ppp)\sm\HH^3$ with coefficients in a local system $T$ is idenfied with
\[\pq{\ZZ[\FF_\Ppp^2-0]\otimes T}\otimes_{\GL_2(\OOO)}M\chain_2\]
and generated in degree $i$ by 
\[\pq{(\alpha,\beta)\otimes t}\otimes[v_1,v_2],\quad(\alpha,\beta)\in\FF_\Ppp^2-0,\quad t\in T.\]

\subsection{Relating the modular and Bianchi complexes to the Galois Lie coalgebra}

\label{sec:relating}

\subsubsection{Motivic correlators at torsion points and averaged base point Hodge correlators}

Let $E$ be an elliptic curve, $p$ a prime, $\Ppp\subset\End E$ a prime over $p$, and $S=E[\Ppp]$. There is an canonical up to root of unity choice of tangent vector $v_0$ at $0\in E$, given by the Dedekind eta function. Extend it to a translation-invariant vector field on $E$ and take $v_s$ to be its fiber at $s$.

Recall that we packaged the Lie coalgebras $\CL^\vee_{E,S,v_s}$ into a coalebra $\CL^\vee_{E,S}=\bigoplus_s\CL^\vee_{E,S,v_s}$. The $\CL^\vee_{E,S,v_s}$ for different $s$ are canonically isomorphic, so there is a natural diagonal $D\subset\CL^\vee_{E,S}$. The image of $D$ under $\Cor_\Hod$ is the space of \emph{averaged base point correlators}. Equivalently, it is the image of the \emph{averaged base point correlator map} $\Cor_\av=\f{1}{\aq{E[\Ppp]}}\sum_s\Cor_s$. The image of the restriction to the space of symmetric correlators is called the coalgebra of \emph{symmetric averaged base point Hodge correlators} and denoted $\Lie_{\rm sym}^\vee(E,E[\Ppp])$. It is the dual to the quotient of $\Lie_\Hod(E,E[\Ppp])$ induced by the quotient of $\gr^W\pi_1^{\rm nil}(E-E[\Ppp],v_0)$ by the adjoint action of $H_1(E;\ZZ)$ and the translation action of $E[\Ppp]$ on $E$.

\subsubsection{Relaxed modular complexes and Hodge correlators}
\label{sec:modular_motivic}
Suppose that $E$ is an elliptic curve, and $\OOO$ its endomorphism ring, and suppose $\Ppp$ is a prime in $\OOO$.$\Ppp$. 

Let $T_k$ denote the graded right $\GL_k(\OOO)$-module $\Sym\pq{H_1(E;\ZZ)^{\oplus k}}\otimes\QQ$, identified with the algebra of polynomials in the variables $t_1,\ol t_1,\dots,t_k,\ol t_k$, and let $\Gamma_1(\Ppp)\subset\GL_k(\OOO)$ be the stabilizer of the vector $(0,\dots,0,1)\in(\OOO/\Ppp)^k$. We will define a map $\theta$ from the relaxed modular complex with coefficients in $T_k$ to the depth $k$ component of the standard cochain complex of the Lie coalgebra $\gr^D\Lie_{\rm sym}^\vee(E,E[\Ppp])$:
\[\theta:T_k\otimes_{\Gamma_1(\Ppp)}\widetilde M_k\ch\to{\rm CE}\ch\pq{\gr^D\Lie_{\rm sym}^\vee(E,E[p])}_{D=k}.\]

Fix an extended basis $\gen{v_1,\dots,v_k,v_0}$ of $V_k(\OOO)$. Also fix an identification of $\FF_\Ppp$ with $E[\Ppp]$.%
We will abuse notation and identify $\alpha\in\FF_\Ppp$ with $\alpha\in E[p]$. Last, we identify the domain of $\theta$ with 
\[\pq{\ZZ[\FF_\Ppp^k-0]\otimes T_k}\otimes_{\GL_2(\OOO)}M_k\ch.\]

In the degree 1 component, define the map on the level of generating series by
\begin{align}
&\sum_{n_1,n_1',\dots,n_k,n_k'}\pq{(\alpha_1,\dots,\alpha_k)\otimes t_1^{n_1}\ol t_1^{n_1'}\dots t_k^{n_k}\ol t_k^{n_k'}}\otimes\bq{v_1,\dots,v_k}\nonumber
\\
&\mapsto
\f{1}{\aq{E[\Ppp]}}\sum_{s\in E[\Ppp]}\Theta_s\du\pg{\alpha_1,\dots,\alpha_k,-\pq{\alpha_1+\dots+\alpha_k}}{t_1:\dots:t_k:0}.
\label{eqn:theta}
\end{align}
The maps in higher degrees are given by
\begin{align*}
&\pq{\ZZ[\FF_\Ppp-0]\otimes T_k}\otimes_{\GL_2(\OOO)}M_k^n\to\pq{\bigwedge^n\gr^D\Lie_{\rm sym}^\vee(E,E[p])}_{D=k},\\
\sum_{n_1,n_1',\dots,n_k,n_k'}&\pq{(\alpha_1,\dots,\alpha_k)\otimes t_1^{n_1}\ol t_1^{n_1'}\dots t_k^{n_k}\ol t_k^{n_k'}}\otimes\pq{\bq{v_1,\dots,v_{k_1}}\wedge\dots\wedge\bq{v_{k_{n-1}+1},\dots,v_{k_n}}}\\
\mapsto\f{1}{\aq{E[\Ppp]}}\sum_{s\in E[\Ppp]}&\Theta_s\du\pg{\alpha_1,\dots,\alpha_{k_1},-(\alpha_1+\dots+\alpha_{k_1})}{t_1:\dots:t_{k_1}:0}\wedge\dots\wedge\\&\wedge\Theta_s\du\pg{\alpha_{k_{n-1}+1},\dots,\alpha_{k_n},-(\alpha_{k_{n-1}+1}+\dots+\alpha_{k_n})}{t_{k_{n-1}+1}:\dots:t_{k_n}:0}.
\end{align*}
\begin{thm}
The map $\theta$ is a well-defined surjective morphism of complexes of graded $\OOO$-modules.
\label{thm:relating}
\end{thm}
\begin{proof}
The map $\theta$ is a morphism of graded $\GL_k(\OOO)$-modules by construction (recall the $t_i,\ol t_i$ are dual to the cohomology generators $\o,\ob$), and is surjective by construction. We need to verify that the map $\theta$ respects (1) the first shuffle relations, (2) the dihedral symmetry relations, (3) the coproduct. We show the three in order.

The first shuffle relation on the image holds termwise -- for each $s\in E[p]$ -- and is equivalent to the relation on the dual generating series (Lemma~\ref{lma:theta_dihedral}(c)):
\[\sum_{\sigma\in\Sigma_{i,j}}\Theta_s\pg{\beta_{\sigma\inv(1)}:\dots:\beta_{\sigma\inv(k)}:0}{u_{\sigma\inv(1)},u_{\sigma\inv(2)},\dots,u_{\sigma\inv(i+j)},-(u_1+\dots+u_{i+j})},\]
where $t_i=u_1+\dots+u_i$, $\alpha_i=\beta_i-\beta_{i-1}$. This is the first shuffle relation on the generating series $\Theta_s\pg{:}{,}$, which holds a priori.

The images of the dihedral symmetry relations are exactly the relations of Corollary~\ref{cly:thetastar_dihedral}, which hold modulo the correlators of elements that are independent of $s$.

Finally, the map $\theta$ intertwines the coproduct. The general case follows from the degree 1. Set $t_0=0$. By (\ref{eqn:thetastar_coproduct}), we have
\begin{align}
\delta\theta\biggl(\sum_{n_1,n_1',\dots,n_k,n_k'}&\pq{(\alpha_1,\dots,\alpha_k)\otimes t_1^{n_1}\ol t_1^{n_1'}\dots t_k^{n_k}\ol t_k^{n_k'}}\otimes\bq{v_1,\dots,v_k}\biggr)=\nonumber\\
=\f{1}{\aq{E_\Ppp}}\sum_{s\in E[\Ppp]}\biggl(\sum_{\rm cyc}\sum_{i=0}^k&\Theta_s\du\pg{\alpha_1,\dots,\alpha_i,-(\alpha_1+\dots+\alpha_i)}{t_1:\dots:t_i:t_0}\wedge\nonumber\\&\wedge\Theta_s\du\pg{\alpha_{i+1},\dots,\alpha_k,-(\alpha_{i+1}+\dots+\alpha_k)}{t_{i+1}:\dots:t_k:t_0}+\text{\rm lower depth terms}\biggr),\label{eqn:dtheta}
\end{align}
where the lower-depth terms are correlators of elements independent of $s$, and the cyclic sum is over indices modulo $k+1$. On the other hand, we have
\begin{align}
\delta\biggl(\sum_{n_1,n_1',\dots,n_k,n_k'}&\pq{(\alpha_1,\dots,\alpha_k)\otimes t_1^{n_1}\ol t_1^{n_1'}\dots t_k^{n_k}\ol t_k^{n_k'}}\otimes\bq{v_1,\dots,v_k}\biggr)\nonumber\\
&=\biggl(\sum_{n_1,n_1',\dots,n_k,n_k'}\pq{(\alpha_1,\dots,\alpha_k)\otimes t_1^{n_1}\ol t_1^{n_1'}\dots t_k^{n_k}\ol t_k^{n_k'}}\otimes\sum_{\rm cyc}\sum_{i=0}^k-[v_1,\dots,v_i]\wedge[v_{i+1},\dots,v_k]\biggr).\label{eqn:justd}
\end{align}
The cyclic shift in $\GL_k(\OOO)$, which maps $v_0\mapsto v_1\mapsto v_2\mapsto\dots\mapsto v_k\mapsto v_0$, acts by the transpose action on the $t_i$ by $t_i\mapsto t_{i+1}-t_1$ (indices modulo $k+1$; recall $t_{k+1}=t_0=0$).
Thus the image of (\ref{eqn:justd}) under $\theta$ agrees with (\ref{eqn:dtheta}) in each summand of the cyclic sum, except with an additive shift of the arguments. It remains to apply the homogeneity of the $\Theta_s\du$.
\end{proof}

\subsubsection*{Remark} Why do we define the map $\theta$ using averaged base point correlators? It would have been possible to define the map to $\Lie_\Hod^\vee(E,E[\Ppp])$ using the correlators with fixed base point, $\Cor_s$. However, this map would be zero. Indeed, any correlator with base point $s$ vanishes modulo the depth filtration in $\Lie_\Hod^\vee(E,E[\Ppp])$ induced by $\Cor_s$, since any correlator can be written modulo those of lower depth by the change of base point formula (\ref{eqn:change_base_point}). Those correlators of lower depth depend on $s$, so this does not imply the image of $\Cor_\av$ is zero.

On the other hand, the map $\theta$ can be modified, replacing $E[\Ppp]$ by its subgroup of order $p$ (if $\aq{E[\Ppp]}=p^2$). We will use this when we specialize $\theta$ to the nodal projective line.

\subsubsection{Bianchi complexes and Hodge correlators in depth 2}

Let $k=2$ and $E$ one of the CM elliptic curves with endomorphism ring $\OOO=\ZZ[i]$ or $\ZZ[\rho]$. According to Lemma~\ref{lma:modular_bianchi}, there is an isomorhpism $\psi:B\ch\to M_2\ch$ from the Bianchi complex to the modular complex. The relaxed modular complex $\widetilde M_2\ch$ is canonically isomorphic to the modular complex $M_2\ch$, since the second shuffle relations are equivalent to the dihedral symmetry relations. Thus we have a map
\[\theta\circ\psi:T_2\otimes_{\Gamma_1(\Ppp)} B\ch\to{\rm CE}\ch\pq{\gr^D\Lie_{\rm sym}^\vee(E,E[\Ppp])}_{D=2}.\]
The complex of the left side is the chain complex with coefficients in the local system $T_2$ on the orbifold $\Gamma_1(\Ppp)\sm\HH^3$. We arrive at the following important result:

\begin{thm}
Let $E$ be one of the CM elliptic curves $E=\CC/\ZZ[i]$ or $E=\CC/\ZZ[\rho]$. Then \[\theta\circ\psi:{\rm CE}\ch\pq{\gr^D\Lie^\vee_\av(E,E[\Ppp])}_{D=2}\] is a surjective morphism of complexes.
\label{thm:bianchi_to_lie}
\end{thm}

It is tempting to extend Thorem~\ref{thm:relating} to higher depth by showing the map $\theta$ descends to the modular complex $M_k\ch$. This requires showing the second shuffle relations for the averaged base point Hod correlators modulo the depth filtration. The following would follow from Conjecture~\ref{cnj:second_shuffle}:
\begin{conj}
The map $\theta$ descends to a morphism of complexes
    \[\theta:T_k\otimes_{\Gamma_1(\Ppp)}M_k\ch\to{\rm CE}\ch\pq{\gr^D\Lie_{\rm sym}^\vee(E,E[\Ppp])}_{D=k}.\]
\end{conj}

\section{Applications}

\label{sec:application}

\subsection{The weight 4 case: Euler complexes}

\label{sec:euler}

Let us show how the map in Theorem~\ref{thm:bianchi_to_lie} generalizes those constructed by \cite{goncharov-levin,goncharov-euler}. To be consistent with those sources, we use the motivic language in this section, but the same results hold in the Hodge realization as well.

\subsubsection{The elements $\theta_E$}

For torsion points $a,b,c\in E$ with $a+b+c=0$, elements $\theta_E(a,b,c)$ are constructed by \cite{goncharov-euler} as follows. 

For $E$ an elliptic curve over a field $k$, $\Ppp\subset\End E$ a prime over $p$, and $z$ a nonzero $\Ppp$-torsion point on $E$, there are elements $\theta_E(z)$, which are $p$-torsion elements in $\ol k_z\du\otimes\ZZ\bq{\f1p}$, where $k_z$ is the extension generated by the coordinates of $z$. They are identified with weight-2 elements in the mixed Tate Lie coalgebra $\Lie_{\MT/\ol k}^\vee$. The real period of the motive $\theta_E(z)$ is $-\log\aq{\theta_E(z)}$.

The elements $\theta_E(a:b:c)$ lie in the Bloch group of $\ol k$, which is identified with the weight-4 part of $\Lie_{\MT/\ol k}^\vee$. We also use the notation $\theta_E(a,b,c)=\theta_E(a:a+b:a+b+c)$, which is unambiguous when $a+b+c=0$ because the $\theta_E(a:b:c)$ are invariant under translation. They are characterized by the following properties:
\begin{enumerate}[(1)]
\item The coproduct is given by \[\delta\theta_E(a,b,c)=\theta_E(a)\wedge\theta_E(b)+\theta_E(b)\wedge\theta_E(c)+\theta_E(c)\wedge\theta_E(a).\]
\item The real period of $\theta_E(a:b:c)$ is given up to a constant multiple by the averaged Chow dilogarithm (\cite{goncharov-arakelov}). The latter can be rewritten as
\begin{equation}
    \f{1}{p^5}\sum_{x\in E[p]}\int_{E(C)}\log\aq{f_{a,x}}\,d^\CC\log\aq{f_{b,x}}\wedge d^\CC\log\aq{f_{c,x}},
    \label{eqn:theta_period}
\end{equation}
where $f_{a,b}$ is a function on $E$ with $\div f_{a,b}=p(\cq{a}-\cq{b})$.
\end{enumerate}

\subsubsection{The elements $\theta_E$ and motivic correlators}

According to \cite{goncharov-hodge-correlators}, \S10.5.5, for $a,b\in E[p]$, the elements $\theta_E(a-b)$ are equal up to a constant multiple to $\Cor_\av\pq{a,b}$. There is a version for the depth 2 elements.
\begin{lma}
Let $E$ be an elliptic curve over a number field. Then, for $a,b,c\in E[p]\sm\cq0$ with $a+b+c=0$, the elements $\theta_E(a:b:c)$ are equal up to a constant multiple to $\Cor_\av(a,b,c)$.
\label{lma:avg_are_theta_e}
\end{lma}
\begin{proof}
The coproduct formulas for the $\theta_E$ and the $\Cor_\av$ concide (\cite{goncharov-hodge-correlators}, Lemma 10.9). It remains to see the periods are equal. Indeed, we take $f_{a,x}$ such that $\log\aq{f_{a,x}(z)}=pG_x(a,z)$, and likewise for $b$ and $c$. Then the formula (\ref{eqn:theta_period}) is evidently a constant multiple of the Hodge correlator
\[\f{1}{\aq{E[p]}}\sum_{x\in E[p]}\Cor_{\HHH,x}(a,b,c),\] as desired.
\end{proof}

The map constructed by \cite{goncharov-euler}, for $\OOO=\ZZ[i]$ or $\ZZ[\rho]$ is:
\begin{align*}
\theta':\ZZ[\Gamma_1(\Ppp)\sm\GL_2(\OOO)]\otimes_{\Gamma_1(\Ppp)} M\ch_2&\to\Lie_\Mot^\vee(E,E[p]),\\
(\alpha_1,\alpha_2)\otimes[v_1,v_2]&\mapsto\theta_E(\alpha_1,\alpha_2,-(\alpha_1+\alpha_2)),\\
(\alpha_1,\alpha_2)\otimes\pq{[v_1]\wedge[v_2]}&\mapsto\theta_E(\alpha_1)\wedge\theta_E(\alpha_2).
\end{align*}

\begin{thm}
The map $\theta'$ is a constant multiple of the component of $\theta\circ\psi$ corresponding to the constant term of the local system $T$.
\end{thm}
\begin{proof}
After unraveling the definitions, in degree 1, this is exactly Lemma~\ref{lma:avg_are_theta_e}, while in degree 2 it amounts to showing that
\[\Cor_\av(0,a)\wedge\Cor_\av(0,b)=\f{1}{\aq{E[p]}}\sum_s\Cor_s(0,a)\wedge\Cor_s(0,b).\]
Expanding the sums and using that $\Cor_s(x,y)\sim\theta_E(x-y)-\theta_E(x-s)-\theta_E(y-s)$ (where we set $\theta_E(0)=0$), this simplifies to
\[\sum_s\theta_E(s)\wedge\theta_E(b-s)+\sum_s\theta_E(a-s)\wedge\theta_E(s)+\sum_s\theta_E(a-s)\wedge\theta_E(b-s)=0.\]
The three sums are both symmetric and antisymmetric under the involutions $s\mapsto b-s$, $s\mapsto a-s$, and $s\mapsto a+b-s$, respectively, so the sum is 0.
\end{proof}

A slight abuse of notation has taken place: $\theta$ maps to $\gr^D\Lie_{\rm sym}^\vee(E,E[\Ppp])$ and $\theta'$ to $\Lie_\MT^\vee(E,E[p])$. However, there is no discrepancy, as the second shuffle relation in weight 4 and depth 2 holds without the lower-depth correction terms, and so the constant term of $\theta$ can be viewed as a map to $\Lie_{\rm sym}^\vee(E,E[p])$. Precisely:
\begin{lma}
For $E$ any elliptic curve and $a,b\in E[\Ppp]$, 
\begin{align*}
    \Cor_\av(a,b,S_{0,0}*S_{0,0})&=0,\\\Cor_\av(a,S_{0,0}*S_{0,0}*S_{0,0})&=0.
\end{align*}
\end{lma}
\begin{proof}
For the first equality, recall that \[S_{0,0}*S_{0,0}=-\f12(\o\otimes\ob-\ob\otimes\o).\] It is easily verified that the coproduct is 0. For the periods, there are two trees contributing to the integral expansion of the correlator. For the tree where $\o,\ob$ are not incident to a common interior vertex, the terms with $\o\otimes\ob$ and with $\ob\otimes\o$ sum to 0. The other tree contributes a constant multiple of
\begin{align*}
&\sum_{s\in E[\Ppp]}\int_EG_s(z,w)\,d^\CC G_s(w,a)\wedge d^\CC G_s(w,b)\wedge \o(z)\wedge\ob(z)\\
&=\sum_{s\in E[\Ppp]}\int_EG_{\rm Ar}(s-w)\, d^\CC G_s(w,a)\wedge d^\CC G_s(w,b)\\
&=\sum_{s\in E[\Ppp]}\int_EG_{\rm Ar}(s-w)\, d^\CC(G_{\rm Ar}(a-w)-G_{\rm Ar}(s-w))\wedge d^\CC(G_{\rm Ar}(b-w)-G_{\rm Ar}(s-w))&=0.
\end{align*}
This follows from the distribution relations for the function $G_{\rm Ar}$, which state that
\[\sum_{s\in E[\Ppp]-0}G_{\rm Ar}(s)=0.\]

The second equality follows simply from dihedral symmetry.
\end{proof}

\subsection{Degeneration to rational curves: Voronoi complexes and multiple $\zeta$-values}
\label{sec:degen}

In this section we study the behavior of the motivic correlators at the boundary of the moduli space $\MMM_{1,n}'$ of elliptic curves with $n$ marked points and a distinguished tangent vector. The results here are a new case of the specialization theorem for correlators on rational curves (\cite{malkin-shuffle}, \S4), and the definitions and proof are analogous. (There is a similar picture for other boundary strata and for higher-genus curves, which can be regarded as the higher-weight version of the results of R.Wentworth \cite{wentworth} about degeneration of Green's functions. We do not expand this subject here.)

\subsubsection{Setup}

It will be enough for us to consider the top boundary stratum in $\ol\MMM_{1,n}'$ in which the elliptic curve $E$ degenerates to nodal $\PP^1$. On an open subset of this stratum, all marked points remain distinct. Furthermore, we will consider degeneration along the direction $\tau=it$, $t\to\infty$ on the modular curve.

Consider an elliptic curve $E$ over $\o B\to\ol\MMM_{1,n}'$, with an open subset $B\to\MMM_{1,n}'$, whose complement $D=\ol B\sm B$ is a normal crossings divisor. A Hodge correlator on $E$ determines a variation of mixed Hodge structures over $B$, which has a canonical extension along every normal vector to $D$. We will describe this canonical extension in the aforementioned case.

The curve $E_\tau\cong\CC/\pq{\ZZ+\ZZ\tau}$ has canonical coordinate $z$, and the nonsingular locus of nodal $\PP^1$ has coordinate $z$ (with the node at $z=0,\infty$) such that 
\[s_a = \begin{cases}z=a&\tau\neq0\\z=e^{2\pi ia}&\tau=0\end{cases},\]
for $a\in\CC$, is a smooth section over $t\in(0,\infty]$. Also fix the relative 1-forms $\o=\f{1}{\sqrt{\Im \tau}}dz,\ob=\f{1}{\sqrt{\Im\tau}}\ol{dz}$ on $E$, which have limit 0 on $\PP^1$.

Let $v_0$ be a relative tangent vector at $s_0$ and $Z_S\subset\CC\du$, $S=\cq{s_a:a\in Z_S}$. Let $\DDD\subset\CL_{E/B,S,v_0}^\vee$ be the subcoalgebra generated by the sections $s_a$, where all the $s_a$ factors are distinct, and the relative 1-forms $\o,\ob$. Also fix the tangent vector $v_0=\f{\del}{\del z}$ at $1\in\PP^1$. 

Let us define a degeneration map
\[\pi_D:\DDD\to\CLie_{\PP^1,Z_S\cup\cq{0},v_0}^\vee\oplus\CLie_{\PP^1,Z_S\cup\cq{\infty},v_0}^\vee.\]
Let $\DDD_T$ be the subspace of $\DDD$ generated by the $s_a$ and elements $\o\otimes\ob$ and $\ob\otimes\o$, which we call the elements of Tate type. If $x\in\DDD$ is not of Tate type, we set $\pi_D(x)=0$. Otherwise, we set 
\begin{align*}
    \pi_D(s_a)&=e^{2\pi i a},\\
    \pi_D(\o\otimes\ob)=-\pi_D(\ob\otimes\o)&=\f12\pq{(0)+(\infty)},
\end{align*}
extended to preserve the tensor product. 

(One can verify, using straightforward but cumbersome combinatorics, that this map is well-defined, i.e., respects the first shuffle relations. This is not required for the results below, since we only require the composition of $\pi_D$ with the Hodge correlator map, a fortiori well-defined.)

\begin{lma}
The map $\pi_D$ is a morphism of coalgebras.
\end{lma}
\begin{proof}
Each term in the coproduct of a generator not of Tate type clearly has a factor that is not of Tate type, because the coproduct preserves the weight. So it is enough to see the map respects the coproduct on the generators of Tate type, considering only the terms of the coproduct where both generators are of Tate type.

Let us do this for the first component of the map, to $\CL^\vee_{\PP^1,Z_s\cup\cq0,v_0}$; the other is analogous. Let $x$ be a generator in $\DD_T$. The coproduct of $x$ has two types of terms: 
\begin{enumerate}[(1)]
    \item the cuts with vertex at some $s_a$ (the term $\delta_S$;
    \item the cuts that give the terms $\delta_{\rm Cas}$.
\end{enumerate}
The coproduct of $\pi_D(x)$ has two types of terms:
\begin{enumerate}[(1$'$)]
    \item the cuts with vertex at some $a\neq0$;
    \item the cuts with vertex at a 0.
\end{enumerate}
The terms (1) that have both factors of Tate type are in obvious bijection with the (1'): observe that the segment that is cut must have the same number of $\o$ and $\ob$ factors on each side -- see Figure~\ref{fig:degeneration_coproduct}, left. Similarly, the terms (2) that have both factors of Tate type are in bijection with the (2') -- see Figure~\ref{fig:degeneration_coproduct}, right.
\end{proof}

\begin{figure}[ht]
\centering
\includegraphics[width=0.3\textwidth]{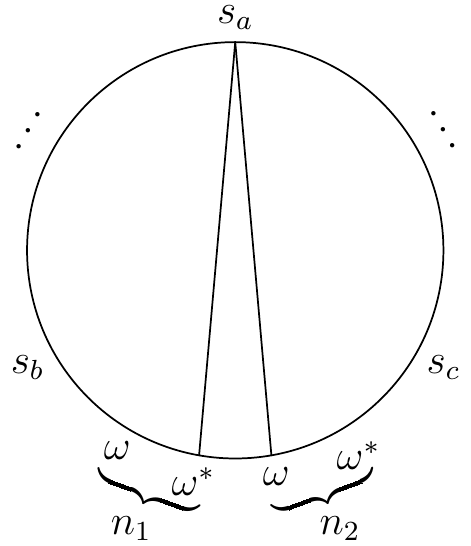}
\hspace{0.2\textwidth}
\includegraphics[width=0.3\textwidth]{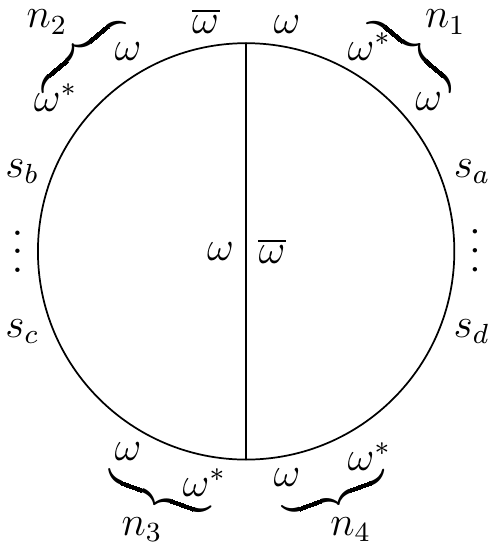}
\\
\includegraphics[width=0.3\textwidth]{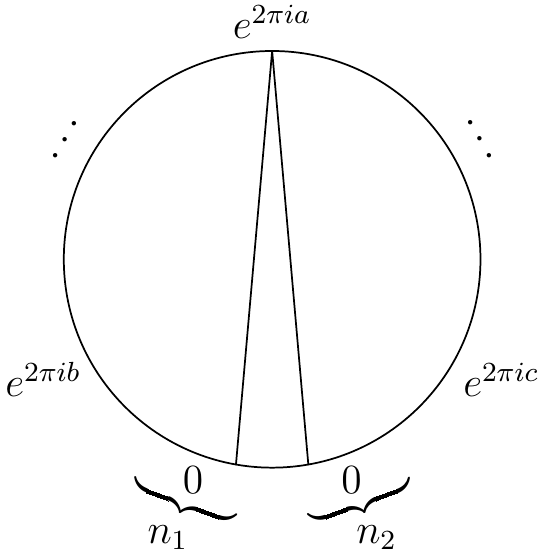}
\hspace{0.2\textwidth}
\includegraphics[width=0.3\textwidth]{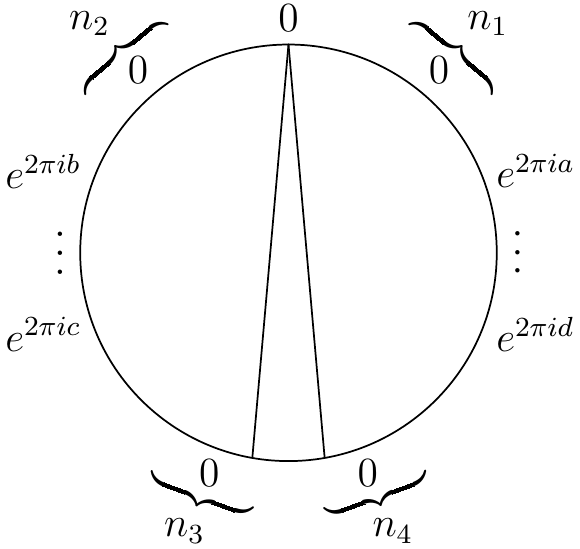}
\simplecap{fig:degeneration_coproduct}{}
\end{figure}

Suppose now that $\ol B,D$ are as above, and that $D$ maps to the boundary stratum in $\ol\MMM_{1,n}'$. Let $\Spec_\infty\Cor_\Hod(x)$ denote the canonical extension of the variation $\Cor_\Hod(x)$ on $B$ to a normal vector to $D$. We then have the following result.
\begin{thm}
    Supose $x\in\DDD$ of weight $n>2$.
    \begin{enumerate}[(a)]
    \item This specialization of the Hodge correlator $\Cor_\Hod(x)$ coincides with the Hodge correlator of the degeneration map:
    \[
      \xymatrix{
        \DDD_{w>2}\ar[r]^{\pi_D\quad\quad\quad\quad\quad\quad}\ar[d]_{\Cor_\Hod}
        &\pq{\CL_{\PP^1,Z_S\cup\cq{0},v_0}^\vee\oplus\CL_{\PP^1,Z_S\cup\cq{\infty},v_0}^\vee}_{w>2}\ar[d]^{\Cor_\Hod}\\
        \Lie_{\Hod/B}^\vee\ar[r]^{\Spec_D}&\Lie_{\Hod}^\vee
      }.
    \]
    \item The Hodge correlator functions on $E$ specialize to the Hodge correlators on $\PP^1$. That is, if $x\in\DDD$ and $\tau=it$, then
    \[
        \lim_{t\to\infty}\Cor_\HHH^{(E_\tau)}(x)\sim\Cor_\HHH^{\PP^1}(\pi_D(x)).    
    \]
    With an appropriate choice of tangent vector on $E_\tau$, this also holds in weight 2.
    \end{enumerate}
\end{thm}
\begin{proof}
We may let $\MMM$ be the moduli space of sets $Z_S$ of $n$ ordered points in $\CC\du$ and $\ol B=(0,\infty]\times\MMM$. We then simultaneously show the following:
\begin{enumerate}[(1)]
\item The periods of $\Cor_\Hod(x)$ extend continuously to $D$.
\item The periods of the specialization of $\Cor_\Hod(x)$ (i.e., the limits of the periods at $D$) coincide with the periods of the degenerations $d(x)$.
\end{enumerate}
The proof is by induction on $w$. Let us see how these imply the result. 

Because the coproduct commutes with specialization, $\Spec_D\Cor_\Hod(x)-\Cor_\Hod(\pi_D(x))$ lies in $\Ext^1_D(\RR(0),\RR(p,q))$, which is one-dimensional and controlled by the period. By (2) it coindices with the period of the degeneration, which immediately gives (1). This implies (a) and (b) in weight $w$.

To show (1), let $q=e^{2\pi i\tau}=e^{-2\pi t}$ be a parameter at the cusp. We show that for $x\in\DDD$, $\Cor_\HHH(x)$ can be represented as a polynomial in $\log q$, where the coefficient of $\log q$ appearing in positive degree has coefficients vanishing at $q=0$ (\emph{tame logarithmic singularities}). This is shown by induction: if $x$ is of weight $w>2$, then $d^\CC\Cor_\HHH(x)$ is expressed in terms of periods of $\delta x$. The latter has logarithmic singularities, by the inductive hypothesis and the fact that the Hodge correlators in weight 1 have logarithmic singularities (see the lemma that follows). Therefore, $\Cor_\HHH(x)$ has tame logarithmic singularities.

By rigidity of $\Ext^1$, we conclude that the difference $\lim_{t\to\infty}\Cor_\HHH^{(E_\tau)}(x)-\Cor_\HHH^{(\PP^1)}(\pi_D(x))$ is independent of the point on $D$, that is, of the choice of $Z_S$. The following lemma implies (3).
\end{proof}

This lemma comprises the analytic ingredients in the preceding proof:
\begin{lma}
\begin{enumerate}[(a)]
    \item Let $x$ be an element of $\DDD$ of weight 2. Then $\Cor_\HHH(x)$ has a logarithmic singularity at $q=0$, and there are constants $c,C$ such that $\lim_{t\to\infty}\Cor_\HHH^{(E_\tau)}(x)-\f Ct-c\Cor_\HHH^{(\PP^1)}(\pi_D(x))=0$.
    \item Let $x=x_0\otimes\dots\otimes x_n$, where each $x_i\in\cq{s_a,\o,\ob}$, be a generator in $\DDD$ of weight $w>2$. Suppose $\lim_{t\to\infty}\Cor_\HHH^{(E_\tau)}(x)-\Cor_\HHH^{(\PP^1)}(\pi_D(x))$ is independent of the choice of $Z_S$. Then this difference is 0.
\end{enumerate}
\end{lma}
\begin{proof}
\begin{enumerate}[(a)]
\item There are three main cases to consider (the remaining ones are symmetric): $x=(s_a)\otimes(s_b)$, $x=(s_a)\otimes\o\otimes\ob$, $x=(s_a)\otimes\o\otimes\o$. The last of those is trivial. For the first two, we use the fact that there is a constant $C$ such that 
\[\lim_{t\to\infty}\pq{G_{\rm Ar}^{(E_{it})}(a)-\f Ct}=\log\aq{\pq{1-e^{2\pi i(a)}}\pq{1-e^{-2\pi i(a)}}}.\]
Therefore, for an appropriate choice of tangent vector $v_t$ at $0\in E_{it}$, we have
\begin{align*}
\lim_{t\to\infty}G_{v_t}^{E_{it}}(a,b)
&=\log\aq{\f{\pq{1-e^{2\pi i(a-b)}}\pq{1-e^{-2\pi i(a-b)}}}{{\pq{1-e^{2\pi ia}}\pq{1-e^{-2\pi ia}}}\pq{1-e^{2\pi ib}}\pq{1-e^{-2\pi ib}}}}\\
&=2\log\aq{\f{e^{2\pi ia}-e^{2\pi ib}}{\pq{1-e^{2\pi ia}}\pq{1-e^{2\pi ib}}}}
&=cG_{v_0}^{(\PP^1)}(a,b),
\end{align*}
Where $v_0=\f{\del}{\del z}$ is a tangent vector at $1\in\PP^1$. This completes the case $x=(s_a)\otimes(s_b)$.

For the case $x=(s_a)\otimes\o\otimes\ob$, notice that $\Cor_\HHH^{(E_{it})}(x)=-G_{\rm Ar}^{(E_{it})}(a)$, so \[\lim_{t\to\infty}\Cor_\HHH^{(E_{it})}(x)=-\log\aq{(1-e^{2\pi ia})(1^{-2\pi ia})}.\]
On the other hand, we also have
\begin{align*}
    G^{(\PP^1)}_{v_0}(z,0)&=\log\aq{\f{z}{1-z}},\\
    G^{(\PP^1)}_{v_0}(z,\infty)&=\log\aq{\f{1}{1-z}},\\
    G^{(\PP^1)}_{v_0}(z,0)+G^{(\PP^1)}_{v_0}(z,\infty)
    &=\log\aq{\f{z}{(1-z)^2}}
    &=-\log\aq{\pq{1-z}\pq{1-1/z}}.
\end{align*}
So we have shown that \[\lim_{t\to\infty}\Cor_\HHH^{(E_{it})}(x)=\f c2\pq{G_{v_0}^{(\PP^1)}(e^{2\pi ia},0)+G_{v_0}^{(\PP^1)}(e^{2\pi ia},\infty)}=c\Cor_\HHH^{(\PP^1)}(\pi_D(x)).\]

\item Let $s_a$ be one of the factors in $x$ (without loss of generality, $x_0=s_a$). We will integrate over $a$ on the segment $[0,1]$. For arbitrary $\tau$, we have
\[\int_{a=0}^1\Cor_\HHH^{(E_\tau)}(x)\,da=0,\] 
since $\int_{a=0}^1 G^{(E_\tau)}(a,z)\,da=0$ for all $z$ by the properties of the Arakelov Green's function, and
\[\int_{a=0}^1\Cor_\HHH^{(\PP^1)}(\pi_D(x))\,da=0,\]
since $\int_{\aq z=1}\Cor_\HHH^{(\PP^1)}(z,b,c)=\LLL_2\pq{\f{z-b}{z-c}}=0$ by the properties of the dilogarithm.
Therefore, \[\int_{a=0}^1\pq{\lim_{t\to\infty}\Cor_\HHH^{(E_\tau)}(x)-\Cor_\HHH^{(\PP^1)}(\pi_D(x))}\,da=0.\] The integrand is independent of $a$, so it is 0.
\end{enumerate}
\end{proof}

\subsubsection{Shuffle relations in depth 2}

Let $x_0,\dots,x_k\in E$ and $m_0,\dots,m_k\geq0$. Define
\[C_{m_0,\dots,m_k}(x_0,\dots,x_k):=\underbrace{S_{0,0}*\dots* S_{0,0}}_{m_0}\otimes(x_0)
\otimes\dots\otimes
\underbrace{S_{0,0}*\dots* S_{0,0}}_{m_k}\otimes(x_k).\]
There is a version of the corrected correlator for this element, where subsets of $\cq{x_0,\dots,x_k}$ are replaced by ``$*$''. We write it in depth 2:
\begin{align*}
    \overline{C}_{m_0,m_1,m_2}(x_0,x_1,x_2):=
    C_{m_0,m_1,m_2}(x_0,x_1,x_2)
    &-C_{m_0,m_1+m_2}(x_0,x_2)
    -C_{m_1,m_2+m_0}(x_1,x_0)
    -C_{m_2,m_0+m_1}(x_2,x_1)\\
    &+C_{m_0+m_1+m_2}(x_0)
    +C_{m_0+m_1+m_2}(x_1)
    +C_{m_0+m_1+m_2}(x_2).
\end{align*}
We have the following variant of Theorem~\ref{thm:dihedral_depth2}:
\begin{lma}
For $m_0,m_1,m_2\geq0$ and $a,b,c\in E$ with $a+b+c=0$,
\[
    \Cor\pq{\ol C_{m_0,m_1,m_2}(a,a+b,a+b+c)+\ol C_{m_0,m_2,m_1}(a,a+c,a+b+c)}=0.
\]
\end{lma}
This is a different version of a second shuffle relation. The proof is identical to that of Theorem~\ref{thm:dihedral_depth2} (we may suppose $a=0$).

Now suppose $a,a+b,a+c,a+b+c$ are distinct and take the correlator with base point 0 over a family with varying $\tau$:
\begin{equation}
    \Cor\pq{\ol C_{m_0,m_1,m_2}(s_a,s_{a+b},s_{a+b+c})+\ol C_{m_0,m_2,m_1}(s_a,s_{a+c},s_{a+b+c})}=0.    
    \label{eqn:variation_of_shuffle}
\end{equation}
An abuse of notation has taken place: the definition of $\ol C$ is assumed to use the relative 1-forms $\o,\ob$ as in the previous section. 

The specialization of the correlator of $\ol C_{m_0,m_1,m_2}(s_a,s_{a+b},s_{a+b+c})$ as $\tau\to i\infty$ is easily seen to be
\begin{align*}
    \Cor(\pi_D(\ol C_{m_0,m_1,m_2}(s_{x_0},s_{x_1},s_{x_2})))
    =&\;
    \Cor_{m_0,m_1,m_2}^{(v_0)}(e^{2\pi ix_0},e^{2\pi ix_1},e^{2\pi ix_2})
    \\&-\f12\biggl(
    \Cor_{m_0,m_1+m_2}^{(v_0)}(e^{2\pi ix_0},e^{2\pi ix_2})
    \\&+\Cor_{m_1,m_2+m_0}^{(v_0)}(e^{2\pi ix_1},e^{2\pi ix_0})
    \\&+\Cor_{m_2,m_0+m_1}^{(v_0)}(e^{2\pi ix_2},e^{2\pi ix_1})\biggr)
    \\&+\Cor^{(v_0)}\text{(terms with $\infty$)}.
\end{align*}
(recall $v_0$ is the tangent vector at $1\in\PP^1$). In particular, by varying $a$ and applying an automorphism of $\PP^1$, we find that the specialization of the relation (\ref{eqn:variation_of_shuffle}) holds for any choice of base point at $\PP^1\sm\cq{0,\infty}$. When it is specialized to $\infty$, the terms with $\infty$ in the specialized correlator vanish. We obtain the relation:
\begin{align*}
    \Cor_{m_0,m_1,m_2}(e^{2\pi ia},e^{2\pi i(a+b)},e^{2\pi i(a+b+c)})
    &-\f12\biggl(
        \Cor_{m_0,m_1+m_2}(e^{2\pi ia},e^{2\pi i(a+b+c)})
        \\&+\Cor_{m_1,m_2+m_0}(e^{2\pi i(a+b)},e^{2\pi ia})
        \\&+\Cor_{m_2,m_0+m_1}(e^{2\pi i(a+b+c)},e^{2\pi i(a+b)})\biggr)\\
    \\+\;\Cor_{m_0,m_1,m_2}(e^{2\pi ia},e^{2\pi i(a+c)},e^{2\pi i(a+b+c)})
    &-\f12\biggl(
        \Cor_{m_0,m_1+m_2}(e^{2\pi ia},e^{2\pi i(a+b+c)})
        \\&+\Cor_{m_1,m_2+m_0}(e^{2\pi i(a+c)},e^{2\pi ia})
        \\&+\Cor_{m_2,m_0+m_1}(e^{2\pi i(a+b+c)},e^{2\pi i(a+c)})
    \biggr)
\end{align*}
correlators now taken with base point at $\infty$. Finally, fix $1=e^{2\pi ia}$ and let $\alpha=e^{2\pi ib}$, $\beta=e^{2\pi ic}$. Rescaling, we arrive at
\[
    \Cor_{m_0,m_1,m_2}(1,\alpha,\alpha\beta)
    +\Cor_{m_0,m_2,m_1}(1,\beta,\alpha\beta)
    -\Cor_{m_0,m_1+m_2}(1,\alpha\beta)
    -\Cor_{m_2+m_0,m_1}(1,\alpha)    
    -\Cor_{m_1+m_0,m_2}(1,\beta).
\]
This is precisely \cite{malkin-shuffle}'s relation (\ref{eqn:p1_second_shuffle}) in depth 2.

\subsubsection{Remark} Let $\mu_p\subset\GG_m$ denote the $p$-th roots of unity. In \cite{goncharov-motivic-modular} (\S2.7), a map from the modular complex for $\GL_2(\ZZ)$ to the standard cochain complex of $\gr^D\Lie_\Hod^\vee(\GG_m,\GG_m-\mu_p)$ is defined using motivic correlators, by a formula similar to (\ref{eqn:theta}):
\[\gamma_2\ch:M_2\ch\otimes_{\Gamma_1(\Ppp)}\QQ\to{\rm CE}\ch\pq{\gr^D\Lie_\Hod^\vee}_2.\]
Alternatively, we can obtain such a map by specializing the map $\theta$. For a generic elliptic curve $E=\CC/(\ZZ+\ZZ\tau)$, we have $\OOO=\ZZ$. We use the variant of the map $\theta$ (\ref{eqn:theta}) defined using an order-$p$ subgroup of $E[\Ppp]$ (see the remark at the end of \S\ref{sec:modular_motivic}). For a family of elliptic curves $E_\tau$ degenerating to nodal $\PP^1$ as $\tau\to+i\infty$, make a continuous choice of such an identification with the subgroup of $E[\Ppp]$ with real coordinates: $z=0,\f1p,\f2p,\dots,\f{p-1}{p}\in E_\tau[p]$. Thus we have a family of sections $s_i$ ($i\in\FF_\Ppp$), which specialize to the $p$-th roots of unity on $\PP^1$. We recover the map $\gamma_2\ch$ by specializing the formula (\ref{eqn:theta}).

\bibliographystyle{alpham}
\bibliography{bibliography-arxiv}

\begin{thebibliography}{GR}

\bibitem[B1]{beilinson-modular}
A.A. Beilinson.
\newblock Higher regulators of modular curves.
\newblock In {\em Applications of algebraic {$K$}-theory to algebraic geometry
  and number theory}, volume~5 of {\em Contemporary Mathematics}, pages 1--34.
  1986.

\bibitem[B2]{beilinson-height-pairing}
A.A. Beilinson.
\newblock Height pairing between algebraic cycles.
\newblock In {\em $K$-theory, arithmetic, and geometry}, volume 1289 of {\em
  LNM}, pages 1--26. Springer, 1987.

\bibitem[B3]{bianchi}
L.~Bianchi.
\newblock Sui gruppi di sostituzioni lineari con coefficienti appartenenti a
  corpi quadratici imaginar{\^i}.
\newblock {\em Mathematische Annalen}, 40:332--412, 1892.

\bibitem[BL]{beilinson-levin}
A.A. Beilinson and A.M. Levin.
\newblock The elliptic polylogarithm.
\newblock In {\em Motives}, volume 55.2 of {\em Symp.\ Pure Math.}, pages
  123--190. AMS, 1994.

\bibitem[D1]{deligne-hodge1}
P.~Deligne.
\newblock {Th\'eorie de Hodge}.
\newblock {\em Actes du congr\`es international des mathématiciens, Nice},
  pages 425--430, 1970.

\bibitem[D2]{deligne-hodge3}
P.~Deligne.
\newblock {Th\'eorie de Hodge: III}.
\newblock {\em Publications Math\'ematiques de l'IH\'ES}, 44:5--77, 1974.

\bibitem[D3]{deligne-p1p1}
P.~Deligne.
\newblock Le groupe fondamental de la droite projective moins trois points.
\newblock In {\em Galois groups over $\QQ$}, volume~16 of {\em Publ. MSRI},
  pages 79--298. MSRI, 1989.

\bibitem[G1]{goncharov-manin}
A.B. Goncharov.
\newblock {The double logarithm and Manin's complex for modular curves}.
\newblock {\em Math. Res. Lett.}, 4:617--636, 1997.

\bibitem[G2]{goncharov-polylogs-modular}
A.B. Goncharov.
\newblock {Multiple polylogarithms, cyclotomy and modular complexes}.
\newblock {\em Math. Res. Lett.}, 5:497--516, 1998.
\newblock \href{http://arxiv.org/abs/1105.2076}{arXiv:1105.2076}.

\bibitem[G3]{goncharov-dihedral}
A.B. Goncharov.
\newblock The dihedral {Lie} algebras and {Galois} symmetries of
  {$\pi^{(l)}_1({\mathbb{P}}^1-(\{0,\infty\}\cup\mu_N))$}.
\newblock {\em Duke Math. Journal}, 110(3):397--487, 2001.
\newblock \href{http://arxiv.org/abs/math/0009121}{arXiv:math/0009121}.

\bibitem[G4]{goncharov-arakelov}
A.B. Goncharov.
\newblock {Polylogarithms, regulators, and Arakelov motivic complexes}.
\newblock {\em Journal of the AMS}, 18:1--60, 2002.
\newblock \href{http://arxiv.org/abs/math/0207036}{arXiv:math/0207036}.

\bibitem[G5]{goncharov-euler}
A.B.\ Goncharov.
\newblock Euler complexes and geometry of modular varieties.
\newblock {\em Geometric and functional analysis}, 17:1872--1914, 2008.
\newblock \href{https://arxiv.org/abs/math/0510310}{arXiv:math/0510310}.

\bibitem[G6]{goncharov-hodge-correlators}
A.B. Goncharov.
\newblock {Hodge correlators}.
\newblock {\em {Journal f\"ur die Reine und Angewandte Mathematik}},
  748:1--138, 2016.
\newblock \href{http://arxiv.org/abs/0803.0297v4}{arxiv:0803.0297v4}.

\bibitem[G7]{goncharov-motivic-modular}
A.B. Goncharov.
\newblock Motivic fundamental group of {${\mathbb{G}}_m\setminus\mu_N$} and
  modular manifolds.
\newblock 2019.
\newblock \href{http://arxiv.org/abs/1910.10321}{arXiv:1910.10321}.

\bibitem[GL]{goncharov-levin}
A.B. Goncharov and A.M. Levin.
\newblock {Zagier's conjecture on $L(E,2)$}.
\newblock {\em Inventiones mathematicae}, 132:393--432, 1995.
\newblock \href{http://arxiv.org/abs/alg-geom/9508008}{arXiv:alg-geom/9508008}.

\bibitem[GR]{goncharov-rudenko}
A.B. Goncharov and D.~Rudenko.
\newblock Motivic correlators, cluster varieties, and {Zagier's} conjecture on
  {$\zeta_{F}(4)$}.
\newblock 2018.
\newblock \href{http://arxiv.org/abs/1803.08585}{arXiv:1803.08585}.

\bibitem[M]{malkin-shuffle}
N.\ Malkin.
\newblock {Shuffle relations for Hodge and motivic correlators}.
\newblock 2020.
\newblock \href{https://arxiv.org/abs/2003.06521}{arXiv:2003.06521}.

\bibitem[W]{wentworth}
R.A. Wentworth.
\newblock {The Asymptotics of the Arakelov-Green's Function and Faltings' Delta
  Invariant}.
\newblock {\em Commun. Math. Phys.}, 137:427--459, 1991.

\end{thebibliography}

\end{document}